\documentclass[a4paper,10pt]{article}

\usepackage{latexsym,amstext,amsmath,amssymb,amsthm,bbm}
\usepackage[all]{xy}

\pdfpagewidth=\paperwidth
\pdfpageheight=\paperheight

\title{Minimal resolutions of geometric $\D$-modules}

\author{R\'emi Arcadias, Universit\'e d'Angers (France)}

\newcommand{\D}{\mathcal{D}}
\newcommand{\C}{\mathbb{C}}

\newcommand{\der}[2]{\frac{\partial #1}{\partial #2}}
\newtheorem{prop}{Proposition}[section]
\newtheorem{theo}{Theorem}[section]
\newtheorem{definition}{Definition}[section]
\newtheorem{lemme}{Lemma}[section]
\newtheorem{cor}{Corollary}[section]

\begin{document}
\maketitle


\begin{abstract}
 
In this paper, we study minimal free resolutions for modules over rings of linear differential operators. The resolutions we are interested in are adapted to a given filtration, in particular to the so-called $V$-filtrations (cf. \cite{res2} and \cite{res1}) .
 We are interested in the module $\D_{x,t}f^s$ associated with germs of functions $f_1,\dots,f_p$, which we call a geometric module, and it is endowed with the $V$-filtration along $t_1=\dots=t_p=0$. The Betti numbers of the minimal resolution associated with this data lead to analytical invariants for the germ of space defined by $f_1,\dots,f_p$ . For $p=1$, we show that under some natural conditions on $f$, the computation of the Betti numbers is reduced to a commutative algebra problem. This includes the case of an isolated quasi homogeneous singularity, for which we give explicitely the Betti numbers. Moreover, for an isolated singularity, we characterize the quasi-homogeneity in terms of the minimal resolution.

\end{abstract}

\section*{Introduction}

In algebraic geometry, minimal free resolutions are a useful tool.
Typically, let $A=\C[x_1,\dots,x_n]$. 
A minimal graded free resolution of a graded module $M$ of finite type is an exact graded sequence
\[
 0\to A^{r_{\delta}}[\mathbf{n}^{(\delta)}]\stackrel{\phi_{\delta}}{\to}\cdots\to A^{r_1}[\mathbf{n}^{(1)}]\stackrel{\phi_1}{\to} A^{r_0}[\mathbf{n}^{(0)}]\stackrel{\phi_0}{\to} M\to 0,
\]
where the entries of the matrices $\phi_i$ for $i\geq 1$ belong to the maximal ideal of $A$. Such a resolution exists and is unique up to graded isomorphisms, so the exponents $r_i$ (called Betti numbers) and the shifts $\mathbf{n}^{(i)}$ are invariants of the graded module $M$.  
In \cite{eisenbud05}, D.Eisenbud gives some geometric applications of these invariants 
 when $M$ is the quotient ring of regular functions of a projective variety. 

Recently, the theory of minimal free resolutions has been studied for algebraic $\D$-modules by T.Oaku et N.Takayama (\cite{res2}). 
Their motivation was the following : $\D$-module theory is an algorithmic tool in algebraic geometry and in linear  partial differential equations. For example, these authors give in \cite{oaku99} 
an algorithm for the computation of the cohomology with complex coefficients of the complement of an hypersurface, using $\D$-modules.
Also, the space of solutions of a system of partial linear differential equations is given by restricting the $\D$-module associated with the system (supposing specialisability), cf. \cite{oaku01}, section 5. 
The point is that in those algorithms, the first thing to do is to compute a free resolution adapted to a filtration.
The authors of \cite{res2} then looked for a way to define free resolutions as short as possible, which led them to the notion of a minimal free resolution.

The minimal free resolutions which we consider are those defined by M.Granger and T.Oaku (\cite{res1}), which work in the analytic local context.

Let $\D$ be the ring of germs of analytic differential operators with analytic coefficients at $0\in\C^n$. Let
\[
(u,v)=(u_1,\dots,u_n,v_1,\dots,v_n)\in\mathbb{Z}^{2n}
\]
such that for any $i,u_i+v_i\geq 0$ and $u_i\leq 0$. Let $x_1,\dots,x_n$ the coordinate functions of $\C^n$. We define a filtration $F^{(u,v)}$ of $\D$ 
by applying the weights $u_i$ to the variables $x_i$ and the weights $v_i$ to the variables $\partial/\partial x_i$.
In particular if for all $i, u_i=0$ and $v_i=1$, this is the classical filtration $F$, denoted $F^{(0,1)}$. 
Take $(u,v)\neq (0,1)$. We define a bifiltration on $\D$ by $F_{d,k}(\D)=F_d^{(0,1)}(\D)\cap F_k^{(u,v)}(\D)$. If $M$ is a $\D$-module 
endowed with a good bifiltration, then $M$ admits bifiltered free resolutions of the type 
\begin{equation}\label{eq0.2}
 \cdots \to\D^{r_2}[\mathbf{n}^{(2)}][\mathbf{m}^{(2)}]\to \D^{r_1}[\mathbf{n}^{(1)}][\mathbf{m}^{(1)}]\to \D^{r_0}[\mathbf{n}^{(0)}][\mathbf{m}^{(0)}]\to M\to 0
\end{equation}
 where $[\mathbf{n}^{(i)}][\mathbf{m}^{(i)}]$ is a vector shift.
We can define, via homogenization, the notion of a minimal bifiltered free resolution. $M$ admits such a resolution, unique up to bifiltered isomorphisms (\cite{res1}, Theorem 3.7). 
Consequently, the Betti numbers $r^i$ and the shifts are invariants of the bifiltered module $M$.

In this article, we study minimal free resolutions of $\D$-modules, adapted to a bifiltration, next we use them to define invariants for singularities of analytic spaces.

In sections 1.1 to 1.4, we recall the construction of minimal bifiltered resolutions, and we focus on a difficulty : to obtain the Betti numbers, we have to compute in a non commutative ring similar to $\D$. The non commutativity makes the effective computation complicated, algorithmically speaking (computations of Gr\"obner bases are heavy) and theoritically speaking because powerful tools of commutative algebra such as Eagon-Nothcott complexes don't always have an equivalent.
In section 1.5, we look for situations in which the minimal bifiltered free resolution of $M$ induces a minimal bigraded free resolution, with same Betti numbers and shifts, of a bigraded module $\textrm{bigr}M$ over the commutative ring
$$\textrm{bigr}(\D)=\bigoplus_{d,k}\frac{F_{d,k}(\D)}{F_{d,k-1}(\D)+F_{d-1,k}(\D)}.$$
Our Theorem \ref{thm1} asserts that we can do this under the following assumption :
\begin{equation}\label{eq0.1}
\forall d,k,\,F_{d,k}(M)=(\cup_{k'} F_{d,k'}(M))\cap (\cup_{d'} F_{d',k}(M)).
\end{equation}
This process will be called the reduction to a commuative algebra problem.

From section 2, we consider the germ of an analytic application
 $f=f_1,\dots,f_p:X\subset (\C^n,0)\to (\C^p,0)$. Let $i_f:X\to X\times\C^p$ be the graph embedding, $\mathcal{O}_X$ the sheaf of analytic functions on $X$, and $\D$ the ring of germs at $0$ of differential operators on $\C^n\times \C^p$. We consider the $\D$-module
$$N_f=(i_{f+}(\mathcal{O}_X))_0$$
which we will call \emph{geometric $\D$-module}. We endow this module with a good bifiltration which refines the $V$-filtration along $\{t_1=\dots=t_p=0\}=X\times 0$, 
we are interested in the Betti numbers and the shifts of the minimal bifiltered free resolution associated with this data.

The module $N_f$ is frequently encountered in the theory of singularities by means of $\D$-modules, and a particular motivation to us is the fact that the algebraic local cohomology $\mathbf{R}\Gamma_{[f=0]}(\mathcal{O}_X)$ is equal to the restriction of $i_{f+}(\mathcal{O}_X)$ along $t_1=\dots=t_p=0$, and the restriction can be obtained from a free resolution adaped to the $V$-filtration (cf. \cite{oaku99}).

Let $\beta_i$ the Betti numbers of $N_f$. Let
 $\beta_f(T)=\sum \beta_i T^i$. Let $n_0$ the embedding dimension at $0$ of the complex space $\mathcal{O}_X/(f_1,\dots,f_p)$, and $r_0$ the cardinal of a minimal set of functions defining this space in $\C^{n_0}$. We show that
\[
 \frac{\beta_f(T)}{(1+T)^{n+p-(n_0+r_0)}}
\]
is a polynomial in $T$ and is an invariant for the complex space $\mathcal{O}_X/(f_1,\dots,f_p)$ (Proposition \ref{prop6}).

From section 2.2, we consider only the hypersurface case, i.e. $p=1$.
We show that for a large class of singularities including those for which the jacobian ideal (generated by $f$ and the derivatives $\partial f/\partial x_i)$ is of linear type (cf. \cite{narvaez}), then the condition (\ref{eq0.1}) 
is satisfied by $N_f$ (Proposition \ref{prop10}).
We also look at the filtered $\D[s]$-module $\D[s]f^s$.
Its Betti numbers are also analytic invariants, and they can be computed by a computer algebra system. 

In section 3, we study quasi homogeneous (for the weights $w_1,\dots,w_n$) isolated singularities. We give explicitely the Betti numbers of $N_f$. Our strategy is, thanks to Theorem \ref{thm1}, to compute in the commutative ring $R=\C\{x\}[t,\tau,\xi_1,\dots,\xi_n]$. The Betti numbers of $N_f$ are the Betti numbers of the quotient $R/I$, where $I$ is the ideal generated by the elements 
\[
 f,t\tau+\sum w_ix_i\xi_i,
\left(\der{f}{x_i}\tau\right)_i,\left(\der{f}{x_i}\xi_j-\der{f}{x_j}\xi_i\right)_{i<j}.
\]
We compute explicitely a minimal bigraded resolution of $R/I$, using Eagon-Northcott complexes.
We obtain the following Betti numbers :
$\beta_{0}=1,\beta_{1}=2+\frac{n(n+1)}{2},\beta_{2}=1+n+2\binom{n+1}{n-2}$, and for $i\geq 3, \beta_{i}=(i-2)\binom{n+2}{i+1}+2\binom{n+1}{i+1}$ (Theorem \ref{thm3}). In particular, let us remark that this sequence of numbers only depends on $n$.

In section 3.2 we characterize the quasi-homogeneity of an isolated singularity by means of the minimal bifiltered resolution of $N_f$. As in commutative algebra, let us define the regularity of a bifiltered module $M$ with respect to the filtration $F$ by $\textrm{reg}_F(M)=\textrm{sup}_{i,j}(\mathbf{n}^{(i)}_j-i)$, considering a minimal bifiltered free resolution of the form (\ref{eq0.2}). We show that $f$ is quasi homogeneous if and only if $\textrm{reg}_F(N_f)=0$ (Proposition \ref{prop20}. In the case of reduced curves, the first Betti numbers is sufficient : if $f$ is quasi homogeneous, then $\beta_1=5$, otherwise $\beta_1\geq 6$ (Proposition \ref{prop5}).

We sincerely thank Michel Granger for his constant help throughout this work, and Luis Narv\'aez-Macarro for many suggestions on the subject. 

\section{Minimal resolutions of $\D$-modules}

\subsection{Rings of differential operators}

The main ring we are interested in is the classical ring of linear differential operators $\D$.
Let $\partial_{x_i}=\partial/\partial x_i$ and for $\beta\in\mathbb{N}^n$, $\partial^{\beta}=\partial_{x_1}^{\beta_1}\dots\partial_{x_n}^{\beta_n}$. 
An element $P$ of $\D$ is written in an unique way
 $P=\sum_{\textrm{finite}}g_{\beta}(x)\partial^{\beta}$
with $g_{\beta}(x)\in \C\{x_1,\dots,x_n\}$. 

A pair $(u,v)=(u_{1},\dots,u_{n},v_{1},\dots,v_{n})\in \mathbb{Z}^{2n}$ is called an 
\emph{admissible weight vector} 
if for all $i$, $\ u_{i}+v_{i}\geq 0$ and $u_{i}\leq 0$. 
If $P=\sum a_{\alpha,\beta}x^{\alpha}\partial^{\beta}\in\D$, let
\[\textrm{ord}_{(u,v)}(P)=\textrm{max}\{\sum u_i\alpha_i+\sum  v_i\beta_i \vert a_{\alpha,\beta}\neq 0\}.\]
We define then a filtration $F_k^{(u,v)}(\D)$ by 
\[F_k^{(u,v)}(\D)=\{P\in\D, \textrm{ord}_{(u,v)}(P)\leq k\}\quad \textrm{for}\ k\in\mathbb{Z}. \] 

\emph{Notation : } For $u=(0,\dots,0)$ and $v=(1,\dots,1)$, we will simply denote by  $F_d(\D)$ the filtration $F_d^{(0,1)}$.

In the sequel we will have to homogenize with respect to the filtration $F$. 
We then consider the graded Rees ring
$\mathbf{R}(\D)=\oplus_d F_d(\D)T^d.$
Let $\D^{(h)}$ be the free $\C$-algebra generated by $\C\{x\}$ and the variables $\partial_{x_1},\dots,\partial_{x_n},h$ quotiented by the relations
\[
 ha-ah,\ h\partial_{x_i}-\partial_{x_i}h,\ \partial_{x_i}\partial_{x_j}-\partial_{x_j}\partial_{x_i},\ \partial_{x_i}a-a\partial_{x_i}-\der{a}{x_i}h
\]
for all $a\in\C\{x\}$. Let us define a graded structure on $\D^{(h)}$
by
\[
 (\D^{(h)})_{d}=\{\sum a_{\beta,k}(x)\partial^{\beta}h^k, \vert\beta\vert + k= d\}
\quad \textrm{for}\, d\in\mathbb{N},\]
and a filtration as for $\D$ by
\[
\textrm{ord}_{(u,v)}(P)=\textrm{max}\{\sum u_i\alpha_i+\sum v_i\beta_i \vert a_{\alpha,\beta,k}\neq 0 \ \textrm{for some}\  k\}
\]
for $P=\sum a_{\alpha,\beta,k}x^{\alpha}\partial^{\beta}h^k\in\D^{(h)}$, i.e. $h$ is given the weight $0$.
We have an isomorphism of graded rings $\D^{(h)}\simeq \mathbf{R}(\D)$
by mapping $a(x)\partial^{\beta}h^k$ to $a(x)\partial^{\beta}T^{k+\vert\beta\vert}$.
%

Let $(u,v)$ be an admissible weight vector. 
We define a bifiltration on $\D$ by
\[
 F_{d,k}(\D)=F^{(0,1)}_d (\D)\cap F^{(u,v)}_k(\D)
\ 
\textrm{for}\, d\in\mathbb{N},k\in\mathbb{Z}
.\]
We will have to deal with the bigraded ring $\textrm{bigr}\D$ defined by
\[
 \textrm{bigr}\D=\bigoplus \textrm{bigr}_{d,k} \D=\bigoplus \frac{F_{d,k}(\D)}{F_{d,k-1}(\D)+F_{d-1,k}(\D)}.
\]

\subsection{Minimal graded (or bigraded) free resolutions}

We recall here the notion of a minimal free resolution for modules over graded or bigraded rings. The rings considered all satisfy a Nakayama lemma, which is a key point in the proof of existence and unicity of minimal resolutions.

\paragraph{Graded case}

Let $A=\bigoplus A_i$ be one of the graded rings $\D^{(h)}$ or 
$\textrm{gr}^F(\D)\simeq\C\{x_1,\dots,x_n\}[\xi_1,\dots,\xi_n]$. 
$A$ possesses an unique graded maximal two-sided ideal $\mathfrak{m}$ generated by $x_1,\dots,x_n$ and $A_1$. Let $M$ be a graded $A$-module of finite type.

For $\mathbf{n}=(n_1,\dots,n_r)\in\mathbb{Z}^r$, we denote by $A^r[\mathbf{n}]$ the graded $A$-module $A^r$ endowed with the graduation $A^r[\mathbf{n}]_d=\bigoplus_{i=1}^r A_{d-n_i}.$
A \emph{graded free resolution}  of $M$ is a graded exact sequence
\begin{equation}\label{eq1.8}
 \cdots \to \mathcal{L}_1\stackrel{\phi_1}{\to}\mathcal{L}_0\stackrel{\phi_0}{\to} M\to 0
\end{equation}
where $\mathcal{L}_i=A^{r_i}[\mathbf{n}^{(i)}]$. 
%
A \emph{minimal} graded free resolution of $M$ is a graded free resolution such that for all $i\geq 1$, $\textrm{im}\phi_i\subset \mathfrak{m}\mathcal{L}_{i-1}$.
Such a resolution exists and is unique up to graded isomorphisms (cf. \cite{eisenbud} Theorem 20.2, or \cite{res1} for non commutative cases). 
Consequently, the exponents $r_i$ and the shifts $\mathbf{n}^{(i)}$ are invariant for the graded module $M$.

\paragraph{Bigraded case}

Here $A=\bigoplus_{d,k}A_{d,k}$ is one of the bigraded rings $\textrm{gr}^V(\D^{(h)})$ or $\textrm{bigr}(\D)$. They possess an unique bigraded two-sided maximal ideal $\mathfrak{m}$. It is generated by $x_1,\dots,x_n$ (where $x_i$ is identified with its class in $A_{0,u_i}$) and $\bigoplus_k A_{(1,k)}$.

Let $r\in\mathbb{N}$ and $\mathbf{n},\mathbf{m}\in \mathbb{Z}^r$. We denote by $A^r[\mathbf{n}][\mathbf{m}]$ the bigraded module $A^r$ with
\[
(A^r[\mathbf{n}][\mathbf{m}])_{d,k}=\bigoplus_{i=1}^r A_{d-n_i,k-m_i}.
\]
A minimal bigraded free resolution of a bigraded module $M$ of finite type is a bigraded exact sequence of the type (\ref{eq1.8}) with $\mathcal{L}_i=A^{r_i}[\mathbf{n}^{(i)}][\mathbf{m}^{(i)}]$ and such that for all $i\geq 1$, $\textrm{im}\phi_i\subset \mathfrak{m}\mathcal{L}_{i-1}$.
Such a resolution exists and is unique up to bigraded isomorphisms, what one can prove in the same way as in the graded case.

\paragraph{Koszul complex}
Classical examples of free resolutions are Koszul complexes.
Let $A$ be a $\C$-algebra, and $a_1,\dots,a_r\in A$ be such that for any $i,j, a_ia_j=a_ja_i$. We define the Koszul complex $K(A;a_1,\dots,a_r)$ :
\[
 0\to A\otimes_{\C} \bigwedge^r \C^r\stackrel{\delta}{\to}
A\otimes_{\C} \bigwedge^{r-1}\C^r\to\cdots\to A\otimes_{\C} \bigwedge^{0}\C^r
\]
with, if $e_1,\dots,e_r$ denotes the canonical basis of $\C^r$ :
\[
\delta(1\otimes e_{i_1}\wedge\dots\wedge e_{i_p})=
\sum_k (-1)^{k-1} a_{i_k}\otimes e_{i_1}\wedge\dots\wedge e_{i_{k-1}}\wedge e_{i_{k+1}}
\wedge\dots\wedge e_{i_p}.
\]
This complex is exact if $A$ is commutative and $a_1,\dots,a_r$ is a regular sequence.
%

\subsection{Filtered free resolutions for $\D$-modules and $\D[s]$-modules}


Let $M$ be a $\D$-module. 
A filtration of $M$ is a collection $\{F_{d}(M)\}_{d\in\mathbb{Z}}$ of $\C\{x\}$-sub modules such that $F_{d}(M)\subset F_{d+1}(M)$, $\bigcup F_{d}(M)=M$, $F_{d'}(\D)F_{d}(M)\subset F_{d+d'}(M)$. 

\begin{definition}
$(F_{d}(M))$ is a \emph{good} filtration if there exist  $f_1,\dots,f_r\in M$ and $n_1,\dots,n_r\in\mathbb{Z}$ such that
\[
 \forall d\in\mathbb{Z}, F_{d}(M)=F_{d-n_1}(\D)f_1+\cdots+F_{d-n_r}(\D)f_r.
\]
\end{definition}
%
For example, if $\mathbf{n}\in \mathbb{Z}^r$, we denote by $\D^r[\mathbf{n}]$ the $\D$-module $\D^r$ endowed with the good filtration 
$F_d(\D^r[\mathbf{n}])=\bigoplus_i F_{d-n_i}(\D)$. 

Suppose that $M$ has a good filtration $(F_d(M))$.
\begin{definition}\label{def1}
 A filtered free resolution of $M$ is a resolution
\begin{equation}\label{eq1.10}
  \cdots \to \mathcal{L}_1\stackrel{\phi_1}{\to}\mathcal{L}_0\stackrel{\phi_0}{\to} M\to 0
\end{equation}
with $\mathcal{L}_i=\D^{r_i}[\mathbf{n}^{(i)}]$, such that for any $d\in\mathbb{Z}$, we have an exact sequence
\[
  \cdots \to F_d(\mathcal{L}_1)\stackrel{\phi_1}{\to}F_d(\mathcal{L}_0)\stackrel{\phi_0}{\to} F_d(M)\to 0.
\]
A minimal filtered free resolution is a filtered free resolution which induces
a minimal graded free resolution of $\textrm{gr}^F (M)$ :
\[
  \cdots \to \textrm{gr}^F\mathcal{L}_1\stackrel{}{\to}\textrm{gr}^F\mathcal{L}_0\stackrel{}{\to} \textrm{gr}^F M\to 0.
\]
\end{definition}

We can also define the notion of a minimal filtered free resolution via homogenization. We have a Rees functor $\mathbf{R}$, from the category of filtered $\D$-modules to the category of graded $\mathbf{R}(\D)$-modules, such that $\mathbf{R}(M)=\oplus F_{d}(M)T^{d}$. We consider a filtered free resolution as (\ref{eq1.10})
which induces a minimal graded free resolution 
\[
  \cdots \to \mathbf{R}\mathcal{L}_1\stackrel{}{\to}\mathbf{R}\mathcal{L}_0\stackrel{}{\to} \mathbf{R}M\to 0.
\]
By \cite{res1} Theorem 3.4, such a resolution exists and is unique up to filtered isomorphisms. 

One easily sees that for a filtered free resolution, it is equivalent to be minimal and to be minimal via homogenization.
Thus minimal filtered free resolution exist and are unique up to filtered isomorphisms.

Now, take a new variable $s$ and consider the ring $\D[s]$. It is endowed with the total filtration
\[
 F_d(\D[s])=\{ \sum P_i s^i, \forall i,\textrm{ord}^F(P_i)+i\leq d\}.
\]
As before, we can define the notion of a good filtration for a $\D[s]$-module, and the notion of a minimal filtered free resolution, which exists and is unique. Here, $\mathfrak{m}$ is the maximal ideal of $\textrm{gr}^F(\D[s])$ generated by $x_1,\dots,x_n$ (identifying $x_i$ with its class in $\textrm{gr}^F_0(\D[s])$), and $\textrm{gr}^F_1(\D[s])$.

\subsection{Bifiltered free resolutions of $\D$-modules}


Let $(u,v)\neq (0,1)$ be an admissible weight vector.

\emph{Notation : } We denote $F_{d}(\D)=F_d^{(0,1)}(\D)$ as before, and $V_{k}(\D)=F_{k}^{(u,v)}(\D)$. 

Let $M$ be a $\D$-module. A \emph{bifiltration} of $M$ is a collection $\{F_{d,k}(M)\}_{d,k\in\mathbb{Z}}$ of $\C\{x\}$-sub modules such that $F_{d,k}(M)\subset F_{d+1,k}(M)\cap F_{d,k+1}(M)$, $\bigcup F_{d,k}(M)=M$, $(F_{d'}(\D)\cap V_{k'}(\D))F_{d,k}(M)\subset F_{d+d',k+k'}(M)$. 

\begin{definition}
$(F_{d,k}(M))$ is a \emph{good} bifiltration if 
there exist  $f_1,\dots,f_l\in M$ and $n_i,m_i\in\mathbb{Z}$ for $i=1,\dots,l$ such that for any $d,k\in\mathbb{Z}$,
\[
 F_{d,k}(M)=(F_{d-n_1}(\D)\cap V_{k-m_1}(\D))f_1+\cdots+(F_{d-n_l}(\D)\cap V_{k-m_l}(\D))f_l.
\]
\end{definition}
For example, if $\mathbf{n},\mathbf{m}\in \mathbb{Z}^r$, 
we denote by $\D^r[\mathbf{n}][\mathbf{m}]$ the $\D$-module $\D^r$ endowed with the good bifiltration 
\[(F_d^{(0,1)}[\mathbf{n}](\D^r)\cap F_k^{(u,v)}[\mathbf{m}](\D^r))_{d,k}.\]
%
Suppose that $M$ is endowed with a good bifiltration. 
Then $M$ has a good $F$-filtration :
$$F_{d}(M)=\bigcup_{k\in\mathbf{Z}} F_{d,k}(M)$$
and a $V$-filtration :
$$V_{k}(M)=\bigcup_{d\in\mathbf{Z}} F_{d,k}(M).$$
$\mathbf{R}(M)$ has a $V$-filtration 
$$V_{k}(\mathbf{R}(M))=\oplus_{d} F_{d,k}(M)T^{d}.$$
The graded $\textrm{gr}^{V}(\D^{(h)})$-module associated with $M$ is
$$\textrm{gr}^{V}(\mathbf{R}(M))=\bigoplus_{k}\textrm{gr}^{V}_{k}(\mathbf{R}(M))\simeq \bigoplus_{k} \bigoplus_{d} \frac{F_{d,k}(M)T^{d}}{F_{d,k-1}(M)T^{d}}.$$

\begin{definition}\label{definition1}
A bifiltered free resolution of $M$ is an exact sequence 
\begin{equation}\label{eq1.11}
\cdots \to \mathcal{L}_{1} \to \mathcal{L}_{0} \to M \to 0
 \end{equation}
with $\mathcal{L}_{i}=\D^{r_{i}}[\mathbf{n}^{(i)}][\mathbf{m}^{(i)}]$, such that 
for any $d,k\in\mathbb{Z}$, we have an exact sequence
\[
 \cdots \to F_{d,k}(\mathcal{L}_{1}) \to F_{d,k}(\mathcal{L}_{0}) \to F_{d,k}(M) \to 0.
\]
A bifiltered free resolution is called minimal if the induced bigraded exact sequence 
$$\cdots \to \textrm{gr}^{V}(\textbf{R}(\mathcal{L}_{1})) \to \textrm{gr}^{V}(\textbf{R}(\mathcal{L}_{0})) \to \textrm{gr}^{V}(\textbf{R}(M)) \to 0 $$
is a minimal bigraded free resolution of $\textrm{gr}^{V}(\textbf{R}(M))$.
\end{definition}
Such a resolution exist and is unique up to bifiltered isomorphism (\cite{res1}, Theorem 3.7). In particular, the exponents $r_{i}$ and the shifts $[\mathbf{n}^{(i)}][\mathbf{m}^{(i)}]$ are invariant for the bifiltered module $M$.

\begin{definition}
 We will call Betti numbers of $(M,(F_{d,k}))$ (resp. shifts) the exponents (resp. shifts) of the minimal bifiltered free resolution of $M$.
\end{definition}

\paragraph{Regularity}
We define the regularity of $M$ with respect to the first filtration.
Consider a minimal bifiltered free resolution of $M$ as (\ref{eq1.11}).
\begin{definition}
$\textrm{reg}_F M=\textrm{sup}_{i,j}(\mathbf{n}^{(i)}_j-i)$.
\end{definition}

We define in the same way the regularity of a bigraded module of finite type over the ring $\textrm{bigr}\D$, with respect to the first graduation.

\subsection{Reduction to a commutative algebra problem}

Let $M$ be a $\D$-module endowed with a good bifiltration. By definition, the Betti numbers of $M$ are those of the bigraded 
 $\textrm{gr}^V(\D^{(h)})$-module $\textrm{gr}^V(\mathbf{R}(M))$. 
Define a bigraded $\textrm{bigr}(\D)$-module
$$\textrm{bigr}(M)=\bigoplus_{d,k} \frac{F_{d,k}(M)}{F_{d-1,k}(N)+F_{d,k-1}(M)}.$$
We will give a condition for $M$ and $\textrm{bigr}M$ to have the same Betti numbers.

Take a bigraded free resolution of $\textrm{gr}^{V}(\textbf{R}(M))$ : 
$$\cdots \to \textrm{gr}^{V}(\textbf{R}(\mathcal{L}_{1})) \to \textrm{gr}^{V}(\textbf{R}(\mathcal{L}_{0})) \to \textrm{gr}^{V}(\textbf{R}(M)) \to 0 $$
with $\mathcal{L}_i=\D[\mathbf{n}^{(i)}][\mathbf{m}^{(i)}]$.
We would like to dehomogenize.
There is a canonical isomorphism of graded rings
$$\textrm{gr}^{V}(\D)\simeq \frac{\textrm{gr}^{V}(\D^{(h)})}{(h-1)}$$
which makes $\textrm{gr}^{V}(\D)$ a left and right $\textrm{gr}^{V}(\D^{(h)})$-module.
We define a dehomogenizing functor $\rho$, from the category of 
bigraded $\textrm{gr}^{V}(\D^{(h)})$-modules to the category of $F$-filtered graded $\textrm{gr}^{V}(\D)$-modules : if $N$ is a bigraded $\textrm{gr}^{V}(\D^{(h)})$-module,  
$$
\rho(N)=\textrm{gr}^V\bigotimes_{\textrm{gr}^V(\D^{(h)})} N\simeq
\frac{N}{(h-1)N},
$$
with the graduation 
\[
 \rho(N)_k =\frac{\bigoplus_d N_{d,k}}{(h-1)\bigoplus_d N_{d,k}}
\]
and the filtration 
\[
 F_d(\rho(N))=1\bigotimes (\bigoplus_k N_{d,k}):=\{1\otimes \omega,
\omega\in \bigoplus_k N_{d,k}\}.
\]
Conversely, we use the canonical isomorphism of rings 
$\mathbf{R}(\textrm{gr}^{V}(\D))\simeq \textrm{gr}^{V}(\D^{(h)})$
to define a homogenizing functor denoted $\mathbf{R}$ as before,
from the category of graded $F$-filtered $\textrm{gr}^{V}(\D)$-modules to the category of bigraded $\textrm{gr}^{V}(\D^{(h)})$-modules.

\begin{definition}
A bigraded $\textrm{gr}^V(\D^{(h)})$-module $N$ is said to be \emph{h-saturated} if $h:N\to N$ is one-to-one. 
\end{definition}

\begin{prop}\label{prop15}
 $\rho$ and $\mathbf{R}$ are exact functors and define an equivalence between the category of bigraded $h$-saturated $\textrm{gr}^{V}(\D^{(h)})$-modules of finite type and the category of
graded $\textrm{gr}^{V}(\D)$-modules endowed with a good $F$-filtration.
\end{prop}
The proof is the same as for Proposition 3.3 of \cite{res1}.

Let us dehomogenize the $\textrm{gr}^V(\D^{(h)})$-module $\textrm{gr}^V(\mathbf{R}(M))$.

\begin{prop}\label{prop0}
There exists an isomorphism of filtered graded $\textrm{gr}^{V}(\D)$-modules 
$$\psi : \textrm{gr}^{V}(M)\simeq\textrm{gr}^{V}(\D)\bigotimes_{\textrm{gr}^{V}(\D^{(h)})} \textrm{gr}^{V}(\textbf{R}(M))$$
where $\textrm{gr}^{V}(M)$ is endowed with the filtration
\[
F_{d}(\textrm{gr}^{V}(M))=
\bigoplus_{k} \frac{F_{d,k}(M)+V_{k-1}(M)}{V_{k-1}(M)}.
\]
\end{prop}
\paragraph{Proof}
Let us make explicit the morphism
\[\psi:\bigoplus\textrm{gr}^V_k(M)\to \textrm{gr}^V(\D)\bigotimes \bigoplus_{d,k} \frac{F_{d,k}(M)T^d}{F_{d,k-1}(M)T^d}.\]
For $\bar{m}\in \textrm{gr}^{V}_{k}(M), m\in F_{d,k}(M)$, we define $\psi(\bar{m})=1\otimes [mT^{d}]$ where $[mT^{d}]$ denotes the class of $mT^d$ in $F_{d,k}(M)T^d/F_{d,k-1}(M)T^d$.
%
The inverse of $\psi$ is the morphism 
\[
 \phi:\textrm{gr}^V(\D)\otimes \textrm{gr}^V(\mathbf{R}(M))\to \textrm{gr}^V(M)
\]
defined as follows : if $P\in V_{k'}(\D)$ and $m\in F_{d,k}(M)$,
$\phi(\overline{P}\otimes[mT^d])=\overline{Pm}\in\textrm{gr}^V_{k+k'}(M)$.

The verifications are straightforward. For example, let us show that $\psi$ is well defined. Let $m\in F_{d_1,k}(M)\cap F_{d_2,k}(M)$ with $d_1<d_2$. Then
\[
 1\otimes [mT^{d_1}]=h^{d_2-d_1}\otimes [mT^{d_1}]=1\otimes [mT^{d_1+d_2-d_1}]=1\otimes [mT^{d_2}]
\]
so $\psi(\bar{m})$ doesn't depend on $d$.
On the other hand, if $m\in V_{k-1}(M)$, then there exists $d$ such that $m\in F_{d,k-1}(M)$ and $[mT^d]=0$ so $\psi (\bar{m})=0$.

Now by definition,
$$F_{d}( \rho({gr}^{V}(\textbf{R}(M)))=1\bigotimes (\textrm{gr}^{V}(\textbf{R}(M)))_{d}=1\bigotimes\bigoplus_{k} \frac{F_{d,k}(M)T^{d}}{F_{d,k-1}(M)T^{d}}.$$
The image of this filtration by $\phi$ defines a filtration on 
 $\textrm{gr}^{V}(M)$ :
$$F_{d}(\textrm{gr}^{V}(M))=\phi (F_{d}(\rho(\textrm{gr}^{V}(\textbf{R}(M))))=
\bigoplus_{k} \frac{F_{d,k}(M)+V_{k-1}(M)}{V_{k-1}(M)}.\square$$

\begin{lemme}\label{lemme13}
The following conditions are equivalent :
\begin{enumerate}
\item $\textrm{gr}^{V}(\textbf{R}(M))$ is $h$-saturated
\item $\forall d,k, F_{d,k}(M)=F_{d+1,k}(M)\cap F_{d,k+1}(M)$
\item $\forall d,k, F_{d,k}(M)=F_{d}(M)\cap V_{k}(M).$
\end{enumerate}
\end{lemme}

\paragraph{Proof}
$(1)\Leftrightarrow(2)$ by definition. $(3)\Rightarrow (2)$ is trivial. 

Let us prove $(2)\Rightarrow (3)$.
Let $x\in F_{d}(M)\cap V_{k}(M)$. We know that there exists $i\geq 0$ such that $x\in F_{d+i,k}(M)\cap F_{d,k+i}(M)$. But condition 2. implies that for all $j\in\mathbb{N}$, $$F_{d,k}=\bigcup_{u+v=j} F_{d+u,k+v}(M)$$ by induction on $j$. Applying this to $j=2i$, we conclude that $x\in F_{d,k}(M)$.
$\square$

Now we can state and prove a result analogous to theorem 3.4 of \cite{res1}, which leads to a definition of minimal bifiltered resolutions via bigraduation.

\begin{theo}\label{thm1}
Let $M$ be a $\D$-module endowed with a good bifiltration such that for any $d,k$, $F_{d,k}(M)=F_{d}(M)\cap V_{k}(M)$. Then $\textrm{bigr}M$ is of finite type over $\textrm{bigr}\D$ and there exists a unique bifiltered free resolution 
$$\cdots \to \mathcal{L}_{1}\to\mathcal{L}_{0}\to M\to 0$$
which induces a minimal bigraded free resolution
$$\cdots \to \textrm{bigr}(\mathcal{L}_{1})\to\textrm{bigr}(\mathcal{L}_{0})\to \textrm{bigr}(M)\to 0.$$
Moreover, this is the minimal bifiltered free resolution of $M$ as defined in Definition \ref{definition1}.
\end{theo}
\paragraph{Proof}
Let us consider the minimal bifiltered free resolution of $M$ :
\begin{equation}
 \label{eq1.4}
\cdots \to \mathcal{L}_{1} \to \mathcal{L}_{0} \to M \to 0
\end{equation}
By definition, this resolution induces a minimal bigraded free resolution
\begin{equation}
 \label{eq1.5}
\cdots \to \textrm{gr}^{V}(\mathbf{R}(\mathcal{L}_{1})) \to \textrm{gr}^{V}(\mathbf{R}(\mathcal{L}_{0})) \to \textrm{gr}^{V}(\mathbf{R}(M)) \to 0.
\end{equation}
Let us dehomogenize (\ref{eq1.5}), i.e. apply 
$\rho$, to obtain a complex
\begin{equation}
 \label{eq1.6}
\cdots \to \mathcal{L}'_{1} \to \mathcal{L}'_{0} \to \textrm{gr}^{V}(M) \to 0
\end{equation}
where $\mathcal{L}'_{i}\simeq \textrm{gr}^{V}(\mathcal{L}_{i})$.

Recall that we assumed $F_{d,k}(M)=F_d(M)\cap V_k(M)$, so $\textrm{gr}^V(\mathbf{R}(M))$ is $h$-saturated by Lemma \ref{lemme13}, this is also true for $\textrm{gr}^V(\mathbf{R}(\mathcal{L}_i))$, so complex (\ref{eq1.6}) is $F$-filtered exact by Proposition \ref{prop15}. By Proposition \ref{prop0}, 
\[F_{d}(\textrm{gr}^{V}(M))=
\bigoplus_{k} \frac{F_{d,k}(M)+V_{k-1}(M)}{V_{k-1}(M)}
\]
so 
\begin{eqnarray*}
\textrm{gr}^{F}_{d}(\textrm{gr}^{V}(M)) & = &
\bigoplus_{k}\frac{F_{d,k}(M)+V_{k-1}(M)}{F_{d-1,k}(M)+V_{k-1}(M)}\\
 & \simeq &
\bigoplus_k \frac{F_{d,k}(M)}{F_{d,k-1}(M)+F_{d-1,k}(M)}\\
 & = & \bigoplus_k \textrm{bigr}_{d,k}M.
\end{eqnarray*}
The same is valid for $\textrm{gr}^V(\mathbf{R}(\mathcal{L}_i))$, then by grading the complex (\ref{eq1.6}), we obtain a bigraded exact sequence
\begin{equation}\label{eq1.7}
\cdots\to\textrm{bigr}(\mathcal{L}_{1})\to\textrm{bigr}(\mathcal{L}_{0})\to\textrm{bigr}(M)\to 0,
\end{equation}
and in particular $\textrm{bigr}M$ is of finite type over $\textrm{bigr}\D$.
To pass from (\ref{eq1.5}) to (\ref{eq1.7}), we substitue $h=0$ in the matrices, so resolution (\ref{eq1.7}) is minimal. We showed that the minimal bifiltered free resolution of $M$ satisfies the theorem. Moreover, $\textrm{bigr}(M)$ and $M$ have the same Betti numbers.

To show the uniqueness, we have to show that a resolution satisfying our theorem is a minimal bifiltered free resolution. Let
\begin{equation}\label{eq1.9}
 \cdots\to\mathcal{L}_1\to \mathcal{L}_0\to M\to 0
\end{equation}
a bifiltered resolution which induces a minimal bigraded resolution 
\[
 \cdots\to \textrm{bigr}(\mathcal{L}_1)\to \textrm{bigr}(\mathcal{L}_0)\to \textrm{bigr}(M)\to 0.
\]
The Betti numbers of resolution (\ref{eq1.9}) are those of $\textrm{gr}^V(\mathbf{R}(M))$ by the first part of the proof.
On the other hand, applying $\mathbf{R}$ to (\ref{eq1.9}) and grading with respect to the $V$-filtration, we obtain a bigraded exact sequence
\[
 \cdots\to \textrm{gr}^V(\mathbf{R}(\mathcal{L}_1))\to \textrm{gr}^V(\mathbf{R}(\mathcal{L}_0))\to \textrm{gr}^V(\mathbf{R}(M))\to 0.
\]
It is thus a minimal bigraded resolution. 
$\square$

\begin{cor}\label{cor3}
Under the assumptions of the previous theorem, the minimal resolutions of the bifiltered module $M$ and the bigradued module $\textrm{bigr}(M)$ have the
same Betti numbers and shifts.
\end{cor}

\section{Minimal resolutions of geometric $\D$-modules}

\subsection{The geometric module $\D_{x,t_1,\dots,t_p}f_1^{s_1}\dots f_p^{s_p}$}

Let $f=(f_1,\dots,f_p):X\subset(\C^n,0)\to (\C^p,0)$ with for any $i, f_i\neq 0$.
Let $\mathcal{O}=\C\{x\}$.
Consider the space 
$$\mathcal{O}\left[\frac{1}{f_1\dots f_p},s_1,\dots,s_p\right]f^s$$
(also denoted $\mathcal{O}[\frac{1}{f},s]f^s$)
where $f^s$ is a symbol representing $f_1^{s_1}\dots f_p^{s_p}$.

Let $\D_{x,t}$ be the ring of differential operators with coefficients in  $\C\{x,t\}$, with $x=(x_1,\dots,x_n)$, $t=(t_1,\dots,t_p)$. We endow $\mathcal{O}[\frac{1}{f},s]f^s$ with a structure of $\D_{x,t}$-module as follows :
the action of $\mathcal{O}$ is trivial ; let $g(x,s)\in \mathcal{O}[\frac{1}{f},s]$, we define
\[
\partial_{x_i}. g(x,s)f^s=\der{g}{x_i}f^s+\sum_j s_j g(x,s)\der{f_j}{x_i}f_j^{-1}f^s,
\]
\[
t_i. g(x,s)f^s=g(x,s_1,\dots,s_i+1,\dots,s_p)f_i f^s,
\]
The action of $t_i$ is bijective, so we define
$\partial_{t_i}=-s_it_i^{-1}$ acting on $\mathcal{O}[\frac{1}{f},s]f^s$.


Our first geometric module is
\[
 N_f=\D_{x,t}f^s=\D_{x,t}f_1^{s_1}\dots f_p^{s_p},
\]
the sub-$\D_{x,t}$-module of $\mathcal{O}[\frac{1}{f},s]f^s$ generated by $f^s$.
Our second geometric module, considered for $p=1$, is $\D[s]f^s$, the sub-$\D[s]$-module of $\mathcal{O}[\frac{1}{f},s]f^s$ generated by $f^s$.

We consider the filtration defined by the weight vector
\[
 (0,\dots,0;-1,\dots,-1;0,\dots,0;1,\dots,1)
\]
corresponding to the variables $(x_1,\dots,x_n;t_1,\dots,t_p;\partial_{x_1},\dots,\partial_{x_n};\partial_{t_1},\dots,\partial_{t_p})$, usually called the $V$-filtration along $t_1=\dots=t_p=0$. We endow $N_f$ with the good bifiltration : 
$F_{d,k}(N_f)=F_{d,k}(\D_{x,t}).f^s$.

We are interested in the Betti numbers and the shifts of the minimal bifiltered resolution of $N_f$, and in their relationship with the singularity defined by $f_1,\dots,f_p$.

\emph{Notation : } Let us denote $(\beta_i(f))_i$ the Betti numbers of the bifiltered module $N_f$. 

Let us give another description of $N_f$.
We denote the $\D$-theoretic direct image of a module $M$ by a morphism $g$ by $g_+(M)$ (cf. \cite{cimpa2}). Consider the graph embedding
\[
i_f : X\to X\times \C^r
\]
defined by $i_f(x)=(x,f_1(x),\dots,f_r(x))$.
We denote by $i_{f+}(\mathcal{O})$ the stalk at $0$ of the sheaf $i_{f+}(\mathcal {O}_X)$. We have an isomorphism of $\C$-vector spaces 

$$i_{f+} (\mathcal{O})\simeq
 \mathcal{O}\otimes_{\C}\C[\partial_{t_1},\dots,\partial_{t_r}]$$
and $i_{f+} (\mathcal{O})$ is generated by $1\otimes 1$ over $\D_{x,t}$. We then take the bifiltration  $F_{d,k}(i_{f+}(\mathcal{O}))=F_{d,k}(\D_{x,t}).(1\otimes 1)$.

Let $I$ be the ideal of $\D_{x,t}$ generated by the elements $(t_j-f_j)_j, (\partial_{x_i}+\sum_j\partial f_j/\partial x_i\partial_{t_j})_i$.
We have isomorphisms
$$\D_{x,t}f^s\simeq i_{f+}(\mathcal{O})\simeq\frac{\D_{x,t}}
{I}$$
with the correspondence $f^s\leftrightarrow 1\otimes 1\leftrightarrow \overline{1}$. 
Indeed, in \cite{malgrange}, Lemme (4.1), it is proved that $I$ is a maximal ideal (written for $p=1$ but the proof works for any $p$).
These isomorphisms are bifiltered by the definition of the bifiltrations.

If one $f_i$ is the zero function, $N_f$ denotes the module $i_{f+}(\mathcal{O})$.

\subsubsection{First properties of the Betti numbers}

\begin{lemme}\label{cor1}
\begin{enumerate}
\item There exists an isomorphism 
$$\mathbf{R}(N_f)\simeq\frac{\D_{x,t}^{(h)}}
{((t_j-f_j)_j, (\partial_{x_i}+\sum_j\der{f_j}{x_i}\partial_{t_j})_i)}$$
which is $V$-filtered if the right-hand-side is endowed with the quotient $V$-filtration.
\item The Betti numbers of the graded $\D_{x,t}^{(h)}$-module $\mathbf{R}(N_f)$ 
are $\binom{n+p}{i}$, for $i=0,\dots,n+p$. (We do not consider here the $V$-filtration)
\end{enumerate}
\end{lemme}

\paragraph{Proof}
1. comes from the bifiltered isomorphism $N_f\simeq \D_{x,t}/I$, and the fact that the elements $(t_j-f_j)_j, (\partial_{x_i}+\sum_j\partial f_j/\partial x_i\partial_{t_j})_i$ form an $F$-involutive basis because their $F$-symbols form a regular sequence in $\textrm{gr}^F(\D_{x,t})$.

Next, the minimal graded resolution of the graded module $\mathbf{R}(N_f)$ is the Koszul complex $K(\D_{x,t}^{(h)};(t_j-f_j)_j, (\partial_{x_i}+\sum_j\partial f_j/\partial x_i\partial_{t_j})_i)$ because the $F^{(0,1)}$-symbols of these elements form a regular sequence in $\textrm{gr}^{(0,1)}(\D_{x,t}^{(h)})$.
$\square$

\begin{prop}
$\forall i, \beta_i(f)\geq \binom{n+p}{i}.$
\end{prop}

\paragraph{Proof}
A minimal bifiltered resolution
\[
\cdots\to\mathcal{L}_1\to\mathcal{L}_0\to N_f\to  0
\]
induces a graded resolution
\[
\cdots\to\mathbf{R}(\mathcal{L}_1)\to\mathbf{R}(\mathcal{L}_0)\to \mathbf{R}(N_f)\to  0
\]
which is a minimal graded $V$-filtered resolution, but not necessarily a minimal graded resolution (forgetting the $V$-filtration). We can minimalize if necessary, which implies the statement by Lemma \ref{cor1}, 2.$\square$

Let us notice that the Betti numbers $\binom{n+p}{i}$ are those corresponding to $(f_1,\dots,f_p)=(x_1,\dots,x_p)$ (one can check it by a standard basis computation). 

\subsubsection{Definition of analytic invariants}

The following results, specially Propositions \ref{prop2} and \ref{prop3}, are inspired by the article \cite{budur} in which the authors define a $b$-function invariant for a complex variety.  

Let $f_1,\dots,f_p\in\mathcal{O}$.
Let $X$ be the germ of space $(\C^n,0)$. 
First, we would like to define invariants for the ideal of $\mathcal{O}$ generated by $f_1,\dots,f_p$.
Suppose $g\in\sum \mathcal{O}f_i$.
Denote $\mathbf{1}\in\mathbb{Z}^{\beta}$ the vector $(1,\dots,1)$.

\begin{prop}\label{prop2}
$\forall i, \beta_i(f,g)=\beta_i(f)+\beta_{i-1}(f)$.
If $[\mathbf{n}^{(i)}][\mathbf{m}^{(i)}]$ is the $i$-th vector shift for $N_f$, 
then the $i$-th vector shift for $N_{f,g}$ is 
$[\mathbf{n}^{(i)},\mathbf{n}^{(i-1)}][\mathbf{m}^{(i)},\mathbf{m}^{(i-1)}]$.
\end{prop}

\paragraph{Proof} 
Let $g=\sum a_{j}f_j$.
Consider the following embeddings :
\begin{description}
\item $i_f : X\to X\times \C^p$
defined by $i_f(x)=(x,f_1(x),\dots,f_p(x))$,
\item $i_{f,g}:X\times \C^{p+1}$
defined by $i_{f,g}(x)=(x,f_1(x),\dots,f_p(x),g(x))$,
\item $i:X\times \C^p\to X\times \C^{p+1}$
defined by $i(x,t)=(x,t,0)$,
\item $\phi:X\times\C^p\to X\times\C^{p+1}$
defined $\phi(x,t)=(x,t,\sum_j a_{j}t_j)$,
\end{description}
and the local isomorphism at the origin  
$\psi:X\times\C^{p+1}\to X\times\C^{p+1}$
defined by $\psi(x,t,u)=(x,t,u+\sum a_{j}t_j)$. 

We have $N_f\simeq  (i_f)_{+} \mathcal{O}$, also
$N_{(f,g)}\simeq (i_{f,g})_{+}\mathcal{O}$. But
$i_{f,g}=\phi\circ i_f=\psi\circ i\circ i_f$,
so 
$N_{(f,g)}\simeq \psi_{+} i_{+} N_f$ by fonctoriality.
$\psi_{+}(\textrm{Id})$ is a bifiltered automorphism of the ring $\D_{x,t,u}$,
 so
$\beta_k(f,g)=\beta_k(i_{+}N_f)$ for any $k$.

Take a minimal bifiltered resolution of $N_f$ :
\begin{equation}\label{eq2.1}
\cdots\to\mathcal{L}_1\to\mathcal{L}_0\to N_f\to 0
\end{equation}
with $\mathcal{L}_i=\D_{x,t}^{l_i}[\mathbf{n}^{(i)}][\mathbf{m}^{(i)}].$

Applying the exact functor $i_{+}$, we get an exact sequence of $\D_{x,t,u}$-modules :
\begin{equation}\label{eq2.2}
\cdots\to i_{+}\mathcal{L}_1\to\ i_{+}\mathcal{L}_0\to
i_{+}N_f\to 0.
\end{equation}
Let us recall that if $M$ is a left $D_{x,t}$-module, then
\[
i_{+}M\simeq M\bigotimes_{\C} \C[\partial_{u}]
\]
with $\D_{x,t,u}$ acting as follows :
\begin{description}
\item $u\bullet m\otimes \partial_u^{\alpha}=-\alpha m\otimes\partial_u^{\alpha-1}$
\item  $\partial_{u}\bullet m\otimes \partial_u^{\alpha}= m\otimes\partial_u^{\alpha+1}$
\item $P(x,\partial_x)\bullet m\otimes \partial_u^{\alpha}=
P(x,\partial_x)m\otimes \partial_u^{\alpha}$.
\end{description}

We have $i_{+}N_f\simeq \psi^{-1}_{+}N_{(f,g)}$ so we get a bifiltration on $i_{+}N_f$ :
\[
F_{d,k}(i_{+}N_f)=F_{d,k}(\D_{x,t})\bullet \delta\otimes 1
=\bigoplus_\xi F_{d-\xi,k-\xi}(\D_{x,t})\delta\otimes \partial_u^{\xi}
\]
where $\delta$ is the canonical generator of $N_f$.
The exact sequence (\ref{eq2.2}) is then bifiltered by defining :
\[
F_{d,k}(\mathcal{L}_i\otimes\C[\partial_u])=
\bigoplus_{\xi}F_{d-\xi,k-\xi}(\mathcal{L}_i)\otimes \partial_u^{\xi}.
\]  
But we have a bifiltered isomorphism
\[
 \mathcal{L}_i\otimes\C[\partial_u]\simeq
\bigoplus_{k=1}^{l_i}\frac{\D_{x,t,u}
[\mathbf{n}_k^{(i)}][\mathbf{m}_k^{(i)}]}{(u)},
\]
so the bifiltration defined is a good one, and if we denote $W$ the ring $\textrm{gr}^V(\D_{x,t,u}^{(h)})$, we have
\[
\textrm{gr}^V \mathbf{R}(\mathcal{L}_i\otimes\C[\partial_u])
\simeq 
\bigoplus_{k=1}^{l_i}\frac{W[\mathbf{n}_k^{(i)}][\mathbf{m}_k^{(i)}]}{(u)}.
\]
Now, let us homogenize and $V$-graduate the exact sequence
(\ref{eq2.2}), we get a $W$-bigraded exact sequence :
\begin{equation}\label{eq2.3} 
\cdots\to \bigoplus_{k=1}^{l_1}\frac{W[\mathbf{n}_k^{(1)}][\mathbf{m}_k^{(1)}]}{(u)}\to
\bigoplus_{k=1}^{l_0}\frac{W[\mathbf{n}_k^{(0)}][\mathbf{m}_k^{(0)}]}{(u)}\to
 \textrm{gr}^V\mathbf{R}(i_{+}N_f)\to 0.
\end{equation}
The minimal bigraded resolution of $W[a][b]/(u)$ is the complex 
\begin{equation}\label{eq2.8}
0\to W[a][b-1]\stackrel{u}{\to} W[a][b]\to 0.
\end{equation}
Seeing (\ref{eq2.3}) as a resolution of $\textrm{gr}^V\mathbf{R}(i_{+}N_f)$ and applying (\ref{eq2.8}) to each of its term, we get a double complex (shifts are omitted) :
\[
\xymatrix{
 W^{l_m} \ar[r] & \cdots\ar[r]& W^{l_1}\ar[r] & W^{l_0}\\
 W^{l_m}\ar[u]\ar[r] & \cdots\ar[r] & W^{l_1}\ar[u]\ar[r] &
 W^{l_0}\ar[u]
}
 \]
The simple complex associated with this double complex gives a bigraded resolution of  $\textrm{gr}^V\mathbf{R}(i_{+}N_f)$. 
Moreover, the entries of the horizontal arrows belong to the maximal ideal because the resolution (\ref{eq2.1}) is supposed to be minimal. The same is true for vertical arrows, then the resolution obtained is minimal, and we get the desired Betti numbers and shifts.
$\square$

Let $\beta_f (T)=\sum_{i\in\mathbb{Z}} \beta_i(f)T^i$.
We interpret Proposition \ref{prop2} as follows :\\
$\beta_{(f,g)}(T)=(1+T)\beta_f(T)$.

\begin{cor}
 The Betti numbers and the shifts of $N_f$ are invariants for the ideal  $I=(f_1,\dots,f_p)$ if $p$ is minimal (i.e. $p=\textrm{dim}_{\C} (I/\mathfrak{m}I)$).
\end{cor}

\paragraph{Proof}
Let $f_1,\dots,f_p$ et $g_1,\dots,g_p$ two minimal generating sets of $I$.
We have $\beta_{(f,g)}(T)=(1+T)^p \beta_f(T)$ and $\beta_{(g,f)}(T)=(1+T)^p \beta_g(T)$. So $\beta_f(T)=\beta_g(T)$. One easily checks that $N_f$ and $N_g$ have the same $i$-th shifts by induction on $i$, by considering the $i$-th shifts of $N_{f,g}$.
$\square$

 We conclude that the Betti numbers and shifts of $N_f$ are invariant for the germs of sub-analytic spaces of $(\C^n,0)$ with $n$ fixed, by associating to such a space its reduced defining ideal.

We would like to define these invariants for germs of complex spaces, indepedently from the embedding. We have to study the dependance of the Betti numbers and the shifts under a closed embedding. Take $y$ a new variable. We will compare $N_f$ with  $N_{(f,y)}$.
Remind that $x=(x_1,\dots,x_n)$ and $t=(t_1,\dots,t_p)$.

\begin{prop}\label{prop3}
 $\forall k, \beta_k(f,y)=\beta_k(f)+2\beta_{k-1}(f)+\beta_{k-2}(f).$
The $i$-th shifts for $N_{(f,y)}$ are
\[
 [n^{(i)},n^{(i-1)},n^{(i-1)}+\mathbf{1},n^{(i-2)}+\mathbf{1}]
[m^{(i)},m^{(i-1)},m^{(i-1)}+\mathbf{1},m^{(i-2)}+\mathbf{1}].
\]
\end{prop}

\paragraph{Proof}
We begin by a technical result using standard bases adapted to the filtration associated with the weights $(u,v)$ (see \cite{ACG} to which we refer for the notions of leading exponent denoted by Exp, and Newton diagram denoted by $\mathcal{N}$).

\begin{lemme}\label{lemme1}
 Let $I$ be an ideal of $\D_{y,u}^{(h)}$ admitting a $V$-adapted standard basis  $P_1,\dots,P_r$ such that for any $i$, $h$ doesn't divide $\textrm{Exp}P_i$, and let $J$ be an ideal of $\D_{x,t}^{(h)}$. Denote $\D^{(h)}$ the ring $\D_{x,y,t,u}^{(h)}$, and $V$ the $V$-filtration restricted to $\D_{x,t}^{(h)}$ and $\D_{y,u}^{(h)}$. Then
\[
 \textrm{gr}^V (\D^{(h)}I+\D^{(h)}J)=\textrm{gr}^V (\D^{(h)})\textrm{gr}^V (I)
+\textrm{gr}^V (\D^{(h)})\textrm{gr}^V (J).
\]
\end{lemme}
\paragraph{Proof}
Here the proof is inspired by the proof of \cite{maisonobe-torelli}, Proposition 4.3.
Let $N$ be the total number of variables $x,t,y,u$.
Take a $V$-adapted standard basis $Q_1,\dots,Q_{r'}$ of $J$. 
Let $U\in \D^{(h)}I+\D^{(h)}J$. It admits a decomposition
\[
 U=\sum U_i P_i+\sum V_j Q_j.
\]
Divide $V_j$ by $P_1, \dots,P_r$ i.e. $V_j=\sum_i U'_{i,j}P_i+W_j.$
So
\[
U=\sum_i (U_i+\sum_j U'_{i,j}Q_j)P_i+\sum W_j Q_j.
\]
We have
\[
 \mathcal{N}(W_j)\bigcap \left(\bigcup_i \textrm{Exp}P_i+\mathbb{N}^{2N+1}\right)=\emptyset.
\]
Also
\[
 \mathcal{N}(W_j Q_j)\bigcap \left(\bigcup_i \textrm{Exp}P_i+\mathbb{N}^{2N+1}\right)=\emptyset
\]
because $Q_j\in\D_{x,t}^{(h)}$ and $h$ doesn't divide $\textrm{Exp}P_i$.
But $P_1,\dots,P_r$ is a standard basis so
\[
 \textrm{Exp}(\sum_i(U_i+\sum_j U'_{i,j}Q_j)P_i)\in(\bigcup_i \textrm{Exp}P_i+\mathbb{N}^{2N+1}).
\]
Then
\[
  \textrm{Exp}(\sum_i(U_i+\sum_j U'_{i,j}Q_j)P_i)\neq
\textrm{Exp}(\sum W_j Q_j).
\]
Thus
\[
 \sigma_V (U)=\epsilon_1 \sigma_V(\sum_i(U_i+\sum_j U'_{i,j}Q_j)P_i)
+\epsilon_2 \sigma_V(\sum W_j Q_j)
\]
with $\epsilon_i=0$ or $1$ for $i=1,2$.$\square$

We have
\[
 \mathbf{R}(N_{(f,y)})=\frac{\D^{(h)}}{\D^{(h)}I+\D^{(h)}J}
\]
where $I$ is the ideal of $\D_{y,u}^{(h)}$ generated by $y-u$ and $\partial_y+\partial_u$, and $J$ is the ideal of $\D_{x,t}^{(h)}$ generated by the elements $(t_j-f_j)_{j=1,\dots,p}$ and 
$(\partial_{x_i}+\sum_j (\partial f/\partial x_j)\partial_{t_j})_{i=1,\dots,n}$.
We can apply the preceding lemma, then if $W$ denotes the ring $\textrm{gr}^V(\D_{x,y,t,u}^{(h)})$, we have 
\begin{eqnarray*}
 \textrm{gr}^V \mathbf{R}N_{(f,y)} & \simeq &
\frac{W}{W \textrm{gr}^V(I)+W \textrm{gr}^V(J)}\\
 & \simeq & 
\frac{W}{W y+W\partial_u+W \textrm{gr}^V(J)}\\
 & \simeq & 
 \textrm{gr}^V\mathbf{R}N_f\otimes_{\C[h]}\C[u,\partial_y,h].
\end{eqnarray*}
Take a minimal bigraded resolution of $\textrm{gr}^V\mathbf{R}N_f$ :
$$\cdots\to\mathcal{L}_1\to\mathcal{L}_0\to \textrm{gr}^V\mathbf{R}N_f\to 0$$
with $\mathcal{L}_i=(\textrm{gr}^V(\D_{x,t}^{h}))^{l_i}[\mathbf{n}^{(i)}][\mathbf{m}^{(i)}]$.\\
Apply the exact functor $-\otimes_{\C[h]} \C[u,\partial_y,h]$. Notice that 
\[
 \mathcal{L}_i\otimes\C[u,\partial_y,h]\simeq 
\bigoplus_{k=1}^{l_i} \frac{W[n_k^{(i)}][m_k^{(i)}]}{(y,\partial_u)},
\]
so we get a $W$-bigraded exact sequence 
\begin{equation}\label{eq2.4}
 \cdots\to\bigoplus_{k=1}^{l_1} \frac{W[n_k^{(1)}][m_k^{(1)}]}{(y,\partial_u)}
\to \bigoplus_{k=1}^{l_0} \frac{W[n_k^{(0)}][m_k^{(0)}]}{(y,\partial_u)}
\to \textrm{gr}^V \mathbf{R}N_{(f,y)}\to 0.
\end{equation}
The minimal bigraded resolution of $W[a][b]/(y,\partial_u)$ is the Koszul complex $K(W;y,\partial_u) : $ 
\[
 0\to W[a+1][b+1]\stackrel{\phi_2}{\rightarrow} W^2[a,a+1][b,b+1] \stackrel{\phi_1}{\rightarrow}W[a][b]\stackrel{\phi_0}{\rightarrow} \frac{W[a][b]}{(y,\partial_u)}\to 0
\]
with $\phi_1(e_1)=y,\phi_1(e_2)=\partial_u$ and
$\phi_2(1)=\partial_u e_1-ye_2$.

Using the exact sequence (\ref{eq2.4}) and the Koxsul complex above, we get a double free complex. The simple complex associated with is the minimal bigraded resolution of $\textrm{gr}^V\mathbf{R}N_{(f,y)}$ and we get the desired Betti numbers and shifts.
$\square$

We interpret Proposition \ref{prop3} as follows :
$\beta_{f,y}(T)=(1+T)^2 \beta_f(T)$.

Let $V$ be a germ of complex space, i.e. a local algebra $\C\{x\}/J$. 
Let $n_0$ be the embedding dimension of $V$ at $0$, and $p_0$ the minimal cardinal of a generating set of the ideal defining $V$ in $\C\{x_1,\dots,x_{n_0}\}$.
Let $f_1,\dots,f_p$ be a generating set of $J$. Let $c_0=p_0+n_0$ (it is an invariant of $V$).

\begin{prop}\label{prop6}
The series
\[
\beta_V(T)=\frac{\beta_f(T)}{(1+T)^{n+p-c_0}}
\]
is a polynomial in $T$ and defines an invariant for the germ of complex space $V$.
\end{prop}

\paragraph{Proof}
The series given is a polynomial because when $n=n_0$ and $p=p_0$, we have $\beta_V(T)=\beta_f(T)$.
To show that this is an invariant, we have to show that $\beta_f(T)/(1+T)^{n+p}$ doesn't depend on the choice of generators of $J$ nor on a closed embedding. It comes from Propositions \ref{prop2} and \ref{prop3}.$\square$

\subsection{The case $p=1$}
 
From now we only consider hypersurfaces, i.e. $p=1$. Let $f\in\C\{x\}$ (different from the zero function). 
The purpose of this section is twofold : first to give a presentation of $\textrm{gr}^V(\mathbf{R}(N_f))$ using $\textrm{ann}_{\D[s]}f^s$, and next 
to give a condition under which we can apply our reduction to commutative algebra to the bifiltered module $N_f$. Finally we make explicit a presentation of $\textrm{bigr}N_f$ in that case.

Let us denote by $f'_i$ the derivative $\partial f/\partial x_i$, for $i=1,\dots,n$.

\subsubsection{Description of $\textrm{gr}^V(\mathbf{R}\D_{x,t}f^s)$ using $\textrm{ann}_{\D[s]}f^s$}

Let $N=\D_{x,t}f^s$ and $W=\textrm{gr}^V(\mathbf{R}(\D_{x,t}))$. We consider $\textrm{gr}^V(\mathbf{R}(N))$ generated by 
$\delta=\overline{f^s T^0}\in \textrm{gr}^V_0(\mathbf{R}(N))$ over $W$.
We want to calculate the annihilator $\textrm{ann}_W \delta$. For $k\geq 0$ we have
\[
 W_{d,k}\simeq F_{d-k}(\D[s])\partial_t^k T^d\quad\textrm{and}\quad
W_{d,-k}\simeq F_{d}(\D[s]) t^k T^d.
\]
Also,
\[
 F_{d,k}(N)=\sum_{i\leq k}F_{d-i}(\D[s])\partial_t^i f^s
\quad\textrm{and}\quad
F_{d,-k}(N)=F_{d}(\D[s])f^{s+k}.
\]
In the following lemma, we determine the operators of the form $P(s)\partial_t^k T^d$ and $P(s)t^k T^d$ which annihilates $\delta$.

\begin{lemme}\label{lemme12}
 \begin{enumerate}
  \item Let $k\geq 1$ and $P(s)\in F_{d-k}(\D[s])$.
Then
$$P(s)\partial_t^k f^s\in\sum_{i\leq k-1} F_{d-i}(\D[s])\partial_t^i f^s\quad \Longleftrightarrow$$
\[
 P(s)\partial_t^k \in\sum_i F_{d-k+1}(\D[s])\partial_t^{k-1}f'_i\partial_t+F_{d-k}(\D[s])\partial_t^k f+\partial_t^k F_{d-k}(\textrm{ann}_{\D[s]}(f^s)).
\]
\item Let $k\geq 0$ and $P(s)\in F_d(\D[s])$.
Then
\[
 P(s)t^k f^s\in F_d(\D[s])f^{s+k+1}\quad \Leftrightarrow\quad
 P(s)t^k\in F_d(\D[s])t^k f+t^k F_{d}(\textrm{ann}_{\D[s]}(f^s)).
\]
\end{enumerate}
\end{lemme}

\paragraph{Proof}
\begin{enumerate}
 \item By induction on $k$. Suppose $k=1$, and $P(s)\partial_t f^s=Q(s)f^s$. We have $\partial_{x_i}f^s=-f'_i\partial_t f^s$, so we can suppose modulo $f'_i\partial_t$ that $Q(s)\in \mathcal{O}[s]$. We have 
\[
 -P(s)sf^{s-1}=Q(s)f^s.
\]
Taking $s=0$, we get $Q(0)(1)=0$ so $Q(0)=0$ because $Q(0)\in\mathcal{O}$. So $Q(s)=s\tilde{Q}(s)$. Simplify by $s$ :
\[
 P(s)f^{s-1}=-\tilde{Q}(s)f^s.
\]
Apply $t$:
\[
 P(s+1)f^s=-\tilde{Q}(s+1)f^{s+1}.
\]
So $P(s+1)=A(s)f+B(s)$ with $B(s)f^s=0$. Multiply on the left by $\partial_t$ :
\[
 P(s)\partial_t=A(s-1)\partial_t f+\partial_t B(s).
\]
Suppose the statement for $k$. Let
\[
 P(s)\partial_t^{k+1}f^s=\sum_{i\leq k}Q_i(s)\partial_t^i f^s.
\]
We can suppose $Q_0(s)\in\mathcal{O}[s]$ by modifying $Q_1(s)$ because of the identity $\partial_{x_i}f^s=f'_i\partial_t f^s$. \\
Take $s=0$. We have 
\[
 \partial_t^i f^s=(-1)^i s(s-1)\cdots (s-i+1)f^{s-i},
\]
then $Q_0(0)(1)=0$ so $Q_0(0)=0$ and $Q_0(s)=s\tilde{Q}_0(s)$. We have 
\[
 -P(s)st^{-1}\partial_t^kf^s=\sum_{1\leq i\leq k}-Q_i(s)st^{-1}\partial_t^{i-1}f^s+s\tilde{Q}_0(s)f^s.
\]
Simplify by $s$, and apply $t$ :
\[
  P(s+1)\partial_t^kf^s=\sum_{1\leq i\leq k}Q_i(s+1)\partial_t^{i-1}f^s-\tilde{Q}_0(s+1)f^{s+1}.
\]
By induction hypothesis,
\[
 P(s+1)\partial_t^k=\partial_t^k A(s)+B(s)\partial_t^k f+\sum C_i(s)\partial_t^k f'_i
\]
with $A(s)f^s=0$.
Multiply to the left by $\partial_t$ :
\[
  P(s)\partial_t^{k+1}=\partial_t^{k+1} A(s)+B(s-1)\partial_t^{k+1} f+\sum C_i(s-1)\partial_t^{k+1} f'_i.
\]
Notice that we conserved the $F$-order during the process.\\
The converse is trivial.

\item Let $P(s)t^kf^s=Q(s)f^{s+k+1}.$ Apply $(t^{-1})^k$ :
\[
 P(s-k)f^s=Q(s-k)f^{s+1}.
\]
Then $P(s-k)=Q(s-k)f+B(s)$ with $B(s)f^s=0$. Multiply by $t^k$ to the left :
\[
 P(s)t^k=Q(s)t^k f+t^k B(s).
\]
The converse is trivial.$\square$
\end{enumerate}

We then have the following :

\begin{prop}\label{prop4}
 $\textrm{ann}_W(\delta)$ is generated by $f,f'_1 \partial_t T,\dots,f'_n \partial_t T$ and $\mathbf{R}(\textrm{ann}_{\D[s]}(f^s))$.
\end{prop}
Here, $\mathbf{R}(\textrm{ann}_{\D[s]}(f^s))$ is an ideal of $\mathbf{R}(\D[s])$ (the Rees ring associated with $\D[s]$ endowed with the filtration $F$). We can see this ideal as a sub-object of $W$ because of the isomorphism $\mathbf{R}(\D[s])\simeq\textrm{gr}^V_0(\mathbf{R}(\D_{x,t})).$



\subsubsection{Condition for the $h$-saturation of $\mathbf{R}(\textrm{gr}^V(N))$}

Let $J(f)$ be the ideal generated by $f'_1,\dots,f'_n$.
Let us define a morphism of graded $\mathcal{O}$-algebras 
\[
\varphi_f : \textrm{gr}^F(\D[s])\simeq\mathcal{O}[s,\xi_1,\dots,\xi_n] \to \bigoplus_{d\geq 0} (\mathcal{O}f+J(f))^dT^d
\]
by $\varphi_f (s)=fT$ and for all $i$, $\varphi_f(\xi_i)=f'_i T$. 
If $P(s)\in\textrm{ann}_{\D[s]}f^s$ with $\textrm{ord}^F(P(s))=d$, one sees that its $F$-symbol $\sigma(P(s))\in\textrm{ker}\,\varphi_f$ by considering the term of degree $d$ in $s$ of $P(s)f^s$ in the space $\mathcal{O}[1/f,s]f^s$. Then in general, $\textrm{gr}^F(\textrm{ann}_{\D[s]}f^s)\subset \textrm{ker}\,\varphi_f$.  

\begin{prop}\label{prop10}
Suppose
\begin{equation}\label{eq2.7}
\textrm{gr}^F(\textrm{ann}_{\D[s]}f^s)=\textrm{ker}\,\varphi_f,
\end{equation}
 then the $W$-module $\textrm{gr}^V(\mathbf{R}(\D_{x,t}f^s))$ is $h$-saturated.
Consequently, the Betti numbers of $N$ are those of the module $\textrm{bigr}N$ by Corollary \ref{cor3}. 
\end{prop}

The condition (\ref{eq2.7}) has been considered in \cite{narvaez09} (property (10)). It holds for a large class of hypersurface singularities, in particular for $f$ such that the ideal $\mathcal{O}f+J(f)$ is of linear type, i.e. for which $\textrm{ker}\,\varphi_f$ is generated by homogeneous elements of degree $1$ (cf. \cite{narvaez09}, Definition 2.2.1). This includes for example quasi homogeneous isolated singularities, which we will study in section 3, and locally quasi homogeneous free divisors (cf. \cite{narvaez}).

\paragraph{Proof}
We have to show : 
\begin{equation}\label{eq2.5}
 \forall d,k, F_{d,k}(N)\cap F_{d+1,k-1}(N)\subset F_{d,k-1}(N).
\end{equation}
For $k\geq 0$ we have
\[
 F_{d,k}(N)=\sum_{l\leq k}F_{d-l}(\D[s])\partial_t^l f^s
\quad\textrm{and}\quad
F_{d,-k}(N)=F_{d}(\D[s])f^{s+k}.
\]
Let us show (\ref{eq2.5}) for $k\leq 0$.
We have to prove : 
if 
$$P(s)f^{l}f^{s}=Q(s)f^{l+1}f^{s}$$
 with $P(s)\in F_{d}(\D[s]), Q(s)\in F_{d+1}(\D[s]), l\in\mathbb{N}$, then $P(s)f^{l}f^{s}=H(s)f^{l+1}f^{s}$, with $H(s)\in F_{d}(\D[s])$.
 We can suppose $l=0$ (apply $t^{-l}$) and $Q\notin F_{d}(\D)$. Applying $t^{-1}$, we get 
$$P(s-1)f^{s-1}=Q(s-1)f^{s}.$$
 Then $\sigma(Q(s-1))\in\textrm{ker}\,\varphi_f$ by considering the terms of degree $d+1$ in $s$ in $\mathcal{O}[1/f,s]f^s$. By the hypothesis, there exists $Q'(s)\in F_{d+1}(\D[s])$ such that $Q'(s)f^s=0$ and $\sigma(Q(s-1))=\sigma(Q'(s))$.
We can take the operator $H(s)=Q(s)-Q'(s+1)$.

Let us prove (\ref{eq2.5}) for $k\geq 1$ by induction on $k$.
Let $k=1$.
Suppose
\[
 P(s)\partial_t f^s=Q(s) f^s
\]
with $P(s)\in F_{d-1}(\D[s])$ and $Q(s)\in F_{d+1}(\D[s])$. 
In the space $\mathcal{O}[1/f,s]f^s$, by considering the terms in $s$ of degree $d+1$ we get $\sigma(Q(s))\in\textrm{ker}\,\varphi_f$ and we conclude like in the case $k=0$, 
making the order of $Q(s)$ smaller.

Assume that (\ref{eq2.5}) is true for $k$ and consider $k+1$.
Suppose
\begin{equation}\label{eq2.6}
 P(s)\partial_t^{k+1}f^s=\sum_{i\leq k}Q_i(s)\partial_t^i f^s
\end{equation}
with $P(s)\in F_{d}(\D[s]), Q_i(s)\in F_{d+1-i}(\D[s])$. 
Using the identity $\partial_{x_i} f^s=-f'_i\partial_t f^s$,
we can assume $Q_0(s)\in\mathcal{O}(s)$ by moving the terms of the form $g(x)\partial^{\beta}s^l$ to $Q_1(s)\partial_t f^s$.
Substitute $s=0$, then $Q_0(0)=0$ and $Q_0(s)=s\tilde{Q}_0(s)$. We have 
\[
 -P(s)st^{-1}\partial_t^kf^s=\sum_{1\leq i\leq k}-Q_i(s)st^{-1}\partial_t^{i-1}f^s+s\tilde{Q}_0(s)f^s.
\]
Simplify by $s$, and apply $t$ :
\[
  P(s+1)\partial_t^kf^s=\sum_{1\leq i\leq k}Q_i(s+1)\partial_t^{i-1}f^s-\tilde{Q}_0(s+1)f^{s+1}.
\]
Then $P(s+1)\partial_t^k f^s\in F_{d-1,k-1}(N)$ by induction. Apply $\partial_t$, so $P(s)\partial_t^{k+1}f^s\in F_{d,k}(N)$.$\square$

\subsubsection{Description of $\textrm{bigr}N$ under $h$-saturation}

\begin{prop}\label{prop11}
If $\textrm{gr}^V(\mathbf{R}(N))$ is $h$-saturated, then the $\textrm{bigr}(\D_{x,t})$-module $\textrm{bigr}(N)$ admits the following presentation :
\[
\textrm{bigr}(N)\simeq\frac{\textrm{bigr}(\D_{x,t})}
{(f,(f'_i\tau),\textrm{gr}^F(\textrm{ann}_{\D[s]}(f^s)))}.
\]
\end{prop}

\paragraph{Proof}
Let $(P_i)_{1\leq i\leq r}$ be an $F$-involutive base of $\textrm{ann}_{\D[s]}f^s$.
By Proposition \ref{prop4}, $N$ admits the following minimal bifiltered presentation (shifts omitted) :
\[
(\D_{x,t})^{1+n+r}\stackrel{\phi_1}{\to}
\D_{x,t}\to N\to 0
\]
where $\phi_1$ sends the canonical basis to the elements $f-t,(f'_i\partial_t+\partial_{x_i}),(P_i)$.
By Theorem \ref{thm1}, this induces the minimal bigraded presentation of $\textrm{bigr}(N)$, the morphism induced by $\phi_1$ mapping the basis elements to $f,(f'_i\tau),\sigma(P_i)$.

\subsubsection{Another geometric module : $\D[s]f^s$}

Now we consider the $\D[s]$-module $\D[s]f^s$ endowed with the good filtration
$F_d(\D[s]f^s)=F_d(\D[s])f^s$. We are also interested in the Betti numbers of this filtered module, and we shall need them in section 3. The advantage is that those numbers are easier to compute with a computer algebra system. On the other hand, we don't know how to extend this invariants to the category of complex spaces (not only hypersurfaces). We suppose $p=1$.

\paragraph{Example : smooth case}
Suppose $f=x_1\in\C\{x_1,\dots,x_n\}$.
One checks that $\textrm{ann}_{\D[s]}(f^s)$ is generated by $x_1\partial_{x_1}-s,\partial_{x_2},\dots,\partial_{x_n}$, 
so
\[
\textrm{gr}^F(\D[s]f^s)\simeq \frac{\mathcal{O}[\xi,s]}{(x_1\xi_1-s,\xi_2,\dots,\xi_n)}.
\]
The minimal graded resolution of $\textrm{gr}^F(\D[s]f^s)$ is the Koszul complex 
\\$K(\mathcal{O}[\xi,s];x_1\xi_1-s,\xi_2,\dots,\xi_n)$, then the Betti numbers are $\binom{n}{i}$, for $i=0,\dots,n$.

Suppose that $f$ is reduced. Take another function $g$. We say that the hypersurfaces $f^{-1}(0)$ and $g^{-1}(0)$ are of same analytic type
if there exists a local analytic isomorphism $\varphi:(\C^n,0)\simeq (\C^n,0)$ such that $\varphi(f^{-1}(0))=g^{-1}(0)$.

\begin{prop}
The Betti numbers and the shifts of the filtered module $\D[s]f^s$ are analytical invariants of $f^{-1}(0)$.
\end{prop}

\paragraph{Proof}
By the Nullstellensatz, it is sufficient to check that the Betti numbers and the shifts are not influenced by a local isomorphism, which is clear, or by multiplying $f$ by an unity, which we shall briefly explain. 
Let $u\in\C\{x\}$ be an unity, let us show that $\D[s]f^s$ and $\D[s](uf)^s$ have the same Betti numbers and shifts. We define a filtered automorphism $\phi$ of the ring $\D[s]$ by 
$\phi_{\vert \mathcal{O}[s]}=\textrm{Id}$ and 
$$\phi(\partial_{x_i})=\partial_{x_i}-su^{-1}\partial_{x_i}(u).$$
One sees that $\phi(\textrm{ann}f^s)=\textrm{ann}(uf)^s$, so $\D[s]f^s$ and $\D[s](uf)^s$ are isomorphic above the ring automorphism $\phi$. $\square$

\section{Isolated singularities and quasi-homogeneity}

Let $f:(\C^{n},0)\to (\C,0)$ be an isolated singularity at the origin.
We will study the Betti numbers of the bifiltered module $N=\D_{x,t}f^{s}$.
 
Denote for  $i=1,\dots,n$, $f'_i=\partial f/\partial x_i$. Let $J(f)$ be the ideal of $\C\{x\}$ generated by $f'_1,\dots,f'_n$, and $\D$ the ring of germs of linear differential operators at the origin in $\C^n$.

An important property of isolated singularities is the fact that $J(f)$ is of linear type. Precisely, consider the morphism of algebras $$\tilde{\varphi}_f:\mathcal{O}[\xi_1,\dots,\xi_n]\to \oplus J(f)^iT^i$$
 which maps $\xi_i$ to $f'_iT$. 

\begin{prop}[\cite{narvaez}, Propositions 3.3 ]\label{prop16}
$\textrm{ker}\,\tilde{\varphi}_f$ is generated by the elements $f'_i\xi_j-f'_j\xi_i$ for $1\leq i<j\leq n$. Consequently, $J(f)$ is of linear type and $\textrm{gr}^F(\textrm{ann}_{\D}f^s)=\textrm{ker}\,\tilde{\varphi}_f$.
\end{prop}

\emph{Notation : }If $1\leq i<j\leq n$, let us denote $S_{i,j}=f'_i\xi_j-f'_j\xi_i$.

\subsection{Isolated quasi homogeneous singularities}

Suppose now that $f$ is quasi homogeneous, i.e. there exist weights $w_{1},\dots,w_{n}$ with $0<w_{i}<1$ for any $i$, such that the Euler vector field 
$\theta=\sum w_{i}x_{i}\partial_{x_{i}}$ satisfies $\theta(f)=f$. 
Denote $\chi=\sum w_i x_i\xi_i$ the symbol of $\theta$.

We are going to compute the Betti numbers $\beta_i(f)$. Our strategy is to use our criterion of reduction to commutative algebra, and to use the filtered $\D[s]$-module  $\D[s]f^s$, for which we also compute the Betti numbers. Let us recall that the morphism of algebras $\phi_f$ maps $s$ to $f$ and $\xi_i$ to $f'_i$.

\begin{prop}\label{prop22}
 $\textrm{ker}\,\varphi_f$ is generated by the elements $s-\chi$ and $f'_i\xi_j-f'_j\xi_i$ for $1\leq i<j\leq n$. Consequently, $\textrm{gr}^F(\textrm{ann}_{\D[s]}f^s)=\textrm{ker}\,\varphi_f$.
\end{prop}

\paragraph{Proof} The first part is a consequence of Proposition \ref{prop16}. 
Then one can show $\textrm{gr}^F(\textrm{ann}_{\D[s]}f^s)=\textrm{ker}\,\varphi_f$ by imitating \cite{narvaez}, Proposition 3.2.$\square$

Then by Proposition \ref{prop10} and Corollary \ref{cor3}, we can work in a commutative ring. 

\emph{Notation : } Let us denote $R=\textrm{bigr}(\D_{x,t})\simeq \C\{x\}[t,\xi,\tau]$ and $\delta$ the class of $f^s\in\textrm{bigr}_{0,0}(N)$.

We have to deal with the $R$-module $\textrm{bigr}(N)$ generated by $\delta$. 
By Propositions \ref{prop11} and \ref{prop22}, we have :

\begin{prop}\label{prop7}
 $\textrm{ann}_R\delta$ is generated by $f,t\tau+\chi, f'_{i}\xi_{j}-f'_{j}\xi_{j} \ \textrm{for}\ 0\leq i<j\leq n$ and $f'_{i}\tau\ \textrm{for}\ 0\leq i\leq n$.
\end{prop}

\subsubsection{Betti numbers of $\D[s]f^s$ and $\D[s]f^s/\D[s]f^{s+1}$ }

Begin with the filtered $\D[s]$-module $\D[s]f^s$. 
By Proposition \ref{prop22}, we have : 
\[
\textrm{gr}^F(\D[s]f^s)=\frac{\textrm{gr}^F(\D[s])}{(s-\chi,(S_{i,j}))}
\simeq \frac{\mathcal{O}[\xi,s]}{(s-\chi,(S_{i,j}))}.
\]
Let us give the minimal graded resolution of the module 
$\mathcal{O}[\xi]/(S_{i,j})_{i<j}.$

For that purpose, let us introduce the \emph{generalized Koszul complexes}, following \cite{northcott}, Appendix C. Let $R$ be a commutative noetherian ring, and $A$ an $m\times n$ matrix with entries in $R$. Denote $\Lambda$ the exterior algebra of the free module $\bigoplus_{i=1}^n R e_i$, $S$ the symmetric algebra of the free module $\bigoplus_{i=1}^m R X_i$ and $\mathfrak{a}$ the ideal generated by all $m$-minors of $A$. For $1\leq i\leq m$, we have a Koszul complex associated with the sequence $a_{i,1},\dots,a_{i,n}$ with differential $\delta_i$. For $1\leq j\leq m$, $X_i$ acts on $S$ by multiplication, and we define a morphism $X_i^{-1}:S\to S$ by dividing the monomials by $X_i$ (the action is zero on a monomial which $X_i$ doesn't divide). Let $t\in\mathbb{Z}$. The generalized Koszul complex $K(A,t)$ is the following : 
\[
 \cdots\to K_h\stackrel{d_h(t)}{\to}K_{h-1}\to\cdots
\]
with
\[
 K_h=\left\{\begin{array}{l}
 \Lambda^{m+h-1}\otimes S_{h-t-1}\quad \textrm{if $h>t$}\\
\Lambda^h\otimes S_{t-h}\quad\textrm{if $h\leq t$}
\end{array}
\right.
\]
and
\[
 d_h(\omega\otimes\alpha)=\left\{\begin{array}{l}
 \sum \delta_i(\omega)\otimes X_i^{-1}(\alpha)\quad \textrm{if $h>t+1$}\\
\delta_m(\delta_{m-1}(\dots \delta_1(\omega)))\otimes\alpha\quad\textrm{if $h=t+1$}\\
 \sum \delta_i(\omega)\otimes X_i(\alpha)\quad \textrm{if $h\leq t$}
\end{array}
\right. .
\]
\begin{theo}[\cite{northcott}, Appendix C, Theorem 2]\label{thm4}
 If $t\leq n-m$ and $0\leq n-m-q+1\leq \textrm{depth}(\mathfrak{a})$, then the truncated complex
\[
 0\to K_{n-m+1}(A,t)\to K_{n-m}(A,t)\to\cdots \to K_q(A,t)
\]
is exact.
\end{theo}

Let us apply this construction to the ring $\mathcal{O}[\xi]$ and the matrix 
\[
 A=\left(\begin{array}{ccc}
 f'_1 & \dots & f'_n\\
-\xi_1 & \dots & -\xi_n
\end{array}\right) .
\]
We have a Koszul complex associated with the regular sequence $f'_1,\dots,f'_n$ with differential $\delta_1$, and another associated with the sequence (also regular) $-\xi_1,\dots,-\xi_n$ with differential $\delta_2$.

The ideal $\mathfrak{a}$ is generated by the $S_{i,j}$. It defines the characteristic variety of the $\D$-module $\D[s]f^s=\D f^s$, which is of dimension at most $n+1$ at all points, by \cite{cimpa1}, Proposition 31. Then the depth of $\mathfrak{a}$ is greater or equal to $2n-(n+1)=n-1$. In Theorem \ref{thm4}, we can take $q\geq n-m+1-(n-1)=0$. We will denote by $K(t)$ the complex $K(A,t)$.

Consider the complex $K(0)$ truncated at degree $0$ :
\[
 \cdots \to\Lambda^4\otimes S_2\to \Lambda^3\otimes S_1\stackrel{d_2(0)}{\to} \Lambda^2\otimes S_0\stackrel{d_1(0)}{\to}\Lambda^0\otimes S_0.
\]
This complex is also called a Eagon-Northcott complex.
It is the minimal graded resolution of $\mathcal{O}[\xi]/(S_{i,j})$ because $d_1(0)(e_i\wedge e_j)=S_{i,j}$. We have $S_i=\bigoplus_{k=0}^i \mathcal{O}[\xi] X_1^k X_2^{i-k}\simeq \mathcal{O}[\xi]^{i+1}$, so the Betti numbers of $\mathcal{O}[\xi]/(S_{i,j})$ are \\ $1,\binom{n}{2},2\binom{n}{3},3\binom{n}{4},\dots,n-1$.

\begin{prop}
The Betti numbers of $\D[s]f^s$ are $\beta_0=1$, $\beta_1=\binom{n}{2}+1$,
and for $i\geq 2, \beta_i=i\binom{n}{i+1}+(i-1)\binom{n}{i}$.
\end{prop}

\paragraph{Proof}
We have
$\textrm{gr}^F(\D[s]f^s)
\simeq \mathcal{O}[\xi,s]/(s-\chi,(S_{i,j}))$
and a short exact sequence
\[
0\to \frac{\mathcal{O}[\xi,s]}{(S_{i,j})}\stackrel{s-\chi}{\rightarrow} \frac{\mathcal{O}[\xi,s]}{(S_{i,j})}\to \frac{\mathcal{O}[\xi,s]}{(s-\chi,(S_{i,j}))}\to 0.
\]
The minimal graded resolution of $\mathcal{O}[\xi,s]/(S_{i,j})$ is the complex $K(0)$. The required minimal graded resolution is then the cone of the morphism of complexes $K(0)\stackrel{s-\chi}{\to}K(0)$. $\square$

Now consider the $\D[s]$-module
\[
 M=\frac{\D[s]f^s}{\D[s]f^{s+1}}
\]
endowed with the good filtration
$F_d(M)=F_d(\D[s])\overline{f^s}.$
We have 
 $M\simeq\textrm{gr}^V_0 (N)$ 
above the ring isomorphism
 $\textrm{gr}^V_0(\D_{x,t})\simeq \D[s].$
Moreover,
$\textrm{gr}^F_d(M)\simeq\textrm{bigr}_{(d,0)}(N).$ 

\begin{prop}\label{prop21}
The Betti numbers of $M$ are
$\beta_0=1,\ \beta_1=2+\binom{n}{2},\ \beta_2= 2\binom{n+1}{3}+1,\ \textrm{and} \ \forall i\geq 3$, 
$\beta_i=(i-2)\binom{n+2}{i+1}+2\binom{n+1}{i+1}$.
\end{prop}

\paragraph{Proof}
By property (\ref{eq2.7}), one can show that we have an exact graded sequence
\[
 0\to \textrm{gr}^F(\D[s]f^s)\stackrel{f}{\to}
 \textrm{gr}^F(\D[s]f^s)\to\textrm{gr}^F(M)\to 0
\]
(see \cite{narvaez09}, Remark 2.2.12).
The minimal graded resolution of $\textrm{gr}^F(M)$ is the cone of the morphism $f$ above, and we know the Betti numbers of $\textrm{gr}^F(\D[s]f^s)$ by the preceding proposition.

\subsubsection{Betti numbers of $N$}

Let $J$ be the ideal of $R$ generated by the elements $f'_{i}\tau$ and $S_{i,j}$.
We have a short exact sequence
$$
0\to \frac{\textrm{ann}\delta}{J}\to\frac{R}{J}\to\frac{R}{\textrm{ann}\delta}\to 0.
$$
We are going to compute the resolutions of $R/J$ and $\textrm{ann}\delta/J$, then we will deduce the Betti numbers of $R/\textrm{ann}\delta$.

\paragraph{Resolution of $R/J$}

Let us begin with the short exact sequence
$$0\to\frac{J}{\sum R S_{i,j}}\to\frac{R}{\sum RS_{i,j}}\to\frac{R}{J}\to 0.$$
We know that the minimal bigraded resolution of 
$R/(S_{i,j})$ is the complex $K(0)$ (with $R$ in place of $\mathcal{O}[\xi]$).
Now let us compute the resolution of $J/\sum R S_{i,j}$. 

\begin{lemme}
 We have an isomorphism
\[
 \frac{\Lambda^1}{\delta_1(\Lambda^2)+\delta_2(\Lambda^2)}\simeq 
\frac{J}{(S_{i,j})}
\]
defined by $\omega\mapsto -\tau\delta_1(\omega)$.
\end{lemme}

\paragraph{Proof} 
Let $\omega\in\Lambda^{1}R^n$ such that $\tau \delta_1(\omega)=\delta_1\circ\delta_2(\eta)$. Taking $\tau=0$, we have 
$\delta_1\circ\delta_2(\eta_{\vert\tau=0})=0$ so $\tau \delta_1(\omega)=\delta_1\circ\delta_2(\eta-\eta_{\vert\tau=0})
\in \tau \delta_1\circ\delta_2( \Lambda^{3})$ then $\delta_1(\omega)=\delta_1\circ\delta_2(\eta')$, so
$\delta_1(\omega-\delta_2(\eta'))=0$.
By the regularity of the sequence $f'_1,\dots,f'_n$, 
the Koszul complex $K(R;f'_1,\dots,f'_n)$ is exact and we have $\omega-\delta_2(\eta')\in \delta_1(\Lambda^2)$.$\square$

Now the minimal resolution of $\Lambda^1/(\delta_1(\Lambda^2)+\delta_2(\Lambda^2))$ is the complex $K(-1)$ truncated at degree $0$ :
\[
 \cdots\to \Lambda^3\otimes S_2\to \Lambda^2\otimes S_1\stackrel{d_1(-1)}{\to}
\Lambda^1\otimes S_0.
\]
We extend the morphism $J/(S_{i,j})\to R/(S_{i,j})$ to a morphism of complexes $\alpha : K(-1)\to K(0)$ defined by $\alpha_0(\omega\otimes 1)=-\tau\delta_1(\omega)\otimes 1$ and for $i\geq 1$, $\alpha_i(\omega\otimes \alpha)=\tau \omega\otimes X_2^{-1}\alpha$. The cone of this morphism is the minimal bigraded resolution of $R/J$.

\paragraph{Resolution of $\textrm{ann}\delta/J$}

\begin{lemme}
We have an isomorphism
\[
\phi_0 : \frac{R^{2}}{((f'_i,-\xi_i),(0,\tau))}\to \frac{\textrm{ann}\delta}{J}
\]
defined by  $X_1\to t\tau+\chi$ and $X_2\to f$.
\end{lemme}

\paragraph{Proof}
We have 
$$-\xi_{k}f+f'_{k}(t\tau+\sum w_{i}x_{i}\xi_{i})=t\tau f'_{k}
+\sum_{i\neq k}w_{i}x_{i}(f'_{k}\xi_{i}-f'_{i}\xi_{k})\in J$$ and
$\tau f-\sum w_{i}x_{i}\tau f'_{i}=0$ so the morphism is well-defined.
Let $M$ be the sub-module of $R^2$ generated by $(0,\tau)$ and the $(f'_i,-\xi_i)$. 
Suppose $(Q,P)\in\textrm{ker}\phi_{0}$ with $P$ and $Q$ bihomogeneous. Necessarily, $\textrm{deg}^{F}(P)\geq 1$ then modulo $M$ we can suppose $P=0$. Then
\begin{equation}\label{eq3.7}
Q.(t\tau+\sum w_{i}x_{i}\xi_{i})=\sum R_{i,j}S_{i,j}+\sum P_{i}f'_{i}\tau.
\end{equation}
Substitue $\tau=0$ and $\xi_i=f'_i$. So $Q_{\vert\tau=0}(\xi_{1}=f'_{1},\dots,\xi_{n}=f'_{n}).f=0$, then $Q_{\vert\tau=0}(\xi_{1}=f'_{1},\dots,\xi_{n}=f'_{n})=0$, and it implies $Q_{\vert\tau=0}\in\sum R S_{i,j}$ by Proposition \ref{prop16}. But $\xi_{i}(f'_j,-\xi_{j})-\xi_{j}(f'_i,-\xi_{i})=(S_{i,j},0)$ so modulo $M$, we can suppose $Q_{\vert\tau=0}=0$.

Taking $\tau=0$ in (\ref{eq3.7}), we get $\sum R_{i,j\vert\tau=0}S_{i,j}=0$ then\\
 $\sum R_{i,j}S_{i,j}\subset \tau\sum R S_{i,j}\subset \sum R f'_i\tau$. So we can suppose $R_{i,j}=0$ by moving $R_{i,j}S_{i,j}$ on $\sum P_{i}f'_{i}\tau$.

Let $Q=Q_{1}\tau+\dots+Q_{d}\tau^{d}$. We have $Q_{d}\tau^{d}t\tau=\sum (P_{i})_{d}f'_{i}\tau \Rightarrow Q_{d}\tau^{d}t=\sum (P_{i})_{d}f'_{i} \Rightarrow 
Q_{d}\tau^{d}=\sum P'_{i}f'_{i}$ because $(f'_{1},\dots,f'_{n},t)$ is a regular sequence. Taking $\tau=0$ we can suppose that $\tau$ divide $P'_{i}$. Then $Q_{d}\tau^{d}\in\sum R f'_{i}\tau$. 
But $\xi_{i}(0,\tau)+\tau(f'_i,-\xi_{i})=(\tau f'_i,0)$ thus modulo $M$, 
we can make the degree in $\tau$ of $Q$ become smaller. We are done when $Q=0$.$\square$

To compute the resolution of $R^2/((f'_i,-\xi_i)_i,(0,\tau))$, we use the following short exact sequence :
$$
0\to\frac{((f'_i,-\xi_i)_i,(0,\tau))}{(f'_i,-\xi_{i})_i}\to
\frac{R^{2}}{(f'_i,-\xi_{i})_i}\to
\frac{R^{2}}{((f'_i,-\xi_i)_i,(0,\tau))}\to 0.
$$
The minimal resolution of $R^{2}/\sum R(f'_i,-\xi_{i})$ is the complex $K(1)$ truncated at degree $0$ (this is valid for $n\geq 3$, but in fact also for $n=2$) :
\[
 \cdots\to \Lambda^4\otimes S_1\to \Lambda^3\otimes S_0\to \Lambda^1\otimes S_0
\stackrel{d_1(1)}{\to} \Lambda^0\otimes S_1\simeq R^2.
\]
This complex is also called a Buchsbaum-Rim complex.

\begin{lemme}
We have an isomorphism 
\[
 \frac{R}{(S_{i,j})_{i<j}}\simeq \frac{((f'_i,-\xi_i)_i,(0,\tau))}{(f'_i,-\xi_{i})_i}
\]
defined by $1\mapsto (0,\tau)$.
\end{lemme}

\paragraph{Proof}
Let $A\in R$ such that $A (0,\tau)=\sum B_i (f'_i,-\xi_i)$. Then $\sum B_i f'_i=0$ so $\sum B_i e_i=\sum_{i<j} R_{i,j} (f'_i e_j-f'_j e_i)$ by regularity of the sequence $f'_1,\dots,f'_n$. So $A\tau=-\sum B_i\xi_i=-\sum R_{i,j} S_{i,j}$ and $A\in \sum R S_{i,j}$.$\square$

We already saw that the minimal resolution of $R/(S_{i,j})_{i<j}$ is the complex $K(0)$, so we have the minimal resolution of the module $((f'_i,-\xi_i)_i,(0,\tau))/(f'_i,-\xi_{i})_i$.
To compute the resolution of $R^2/((f'_i,-\xi_i)_i,(0,\tau))$, 
it is thus sufficient to extend the morphism $((f'_i,-\xi_i)_i,(0,\tau))/(f'_i,-\xi_{i})_i\to
 R^{2}/(f'_i,-\xi_{i})_i$ to a morphism of resolutions, also called $\alpha$. 
We can take $\alpha_0 : \Lambda^0\otimes S_0\to \Lambda^0\otimes S_1$ defined by $\alpha_0(1\otimes 1)=\tau\otimes X_2$, $\alpha_1:\Lambda^2\otimes S_0\to \Lambda^1\otimes S_0$ such that $\alpha_1(\omega\otimes 1)=\tau\delta_1(\omega)\otimes 1$, and for $i\geq 2$, $\alpha_i:\Lambda^{i+1}\otimes S_{i-1}\to \Lambda^{i+1}\otimes S_{i-2}$ defined by $\alpha_i(\omega\otimes\eta)=-\tau\omega\otimes X_2^{-1}\eta$.
The cone of this morphism is the minimal resolution of $R^2/((f'_i,-\xi_i)_i,(0,\tau))\simeq\textrm{ann}\delta/J$.

\paragraph{Betti numbers of $R/\textrm{ann}\delta$}
\begin{theo}\label{thm3}
\begin{enumerate}
 \item 
 The Betti numbers of $N$ are :
$\beta_{0}=1,\ \beta_{1}=2+\frac{n(n+1)}{2},\ \beta_{2}=1+n+2\binom{n+1}{n-2}$,\\ and for $i\geq 3, \beta_{i}=(i-2)\binom{n+2}{i+1}+2\binom{n+1}{i+1}.$
\item $\textrm{reg}_F N=0$.
\end{enumerate}
\end{theo}
For example, for $n=2$ the Betti numbers are $1,5,5,1$ and for $n=3$ they are $1,8,12,7,2$. That was conjectured by the authors of \cite{res1} (see the examples they give).

\paragraph{Proof} Let us prove 1.
We know the minimal resolutions of $\textrm{ann}\delta/J$ and $R/J$. Denote them respectively $L$ and $L'$. We can extend the morphism $\textrm{ann}\delta/J\to R/J$ to a morphism of resolutions $\alpha:L\to L'$. The cone of this morphism is then a resolution of $R/\textrm{ann}\delta$. Unfortunately, it is not easy to write explicitely all the arrows $\alpha_{i}$.
Let us give them for $i=0$ and $i=1$.
$\alpha_{0}:\Lambda^0\otimes S_1\to \Lambda^0\otimes S_0$
is defined by $\alpha_{0}(1\otimes X_1)=(t\tau+\chi)\otimes 1$ and $\alpha_{0}(1\otimes X_2)=f\otimes 1$. Next, we can take
\[
\alpha_1 : \Lambda^0\otimes S_0\oplus \Lambda^1\otimes S_0\to
 \Lambda^1\otimes S_0\oplus \Lambda^2\otimes S_0
\]
defined by 
$\alpha_1(1\otimes 1,0)=(-\sum w_i x_i e_i,0)$ and
$\alpha_1(0,\omega\otimes 1)=(-t\omega,-\omega\wedge \sum w_i x_i e_i)$.

The cone of $\alpha$ is a bigraded resolution of $\textrm{bigr}(N)$ 
\[
 \cdots\to R^{\beta_2}\stackrel{\Phi_2}{\longrightarrow}R^{\beta_1}\stackrel{\Phi_1}{\longrightarrow}R\stackrel{\Phi_0}{\longrightarrow}\textrm{bigr}(N)\to 0
\]
whose Betti numbers are those of the statement of our theorem. 
Let us show that this is the minimal bigraded resolution of $\textrm{bigr}N$, i.e. for all $i\geq 1$, $\textrm{im}\Phi_i\subset\mathfrak{m}R^{\beta_{i-1}}$. For $i=1$ and $i=2$, that comes from the explicit computation of the arrows $\alpha_0$ and $\alpha_1$, whose entries belong to the maximal ideal.

 For $i\geq 3$, we proceed as follows : suppose $\textrm{im}\Phi_{i_0}\nsubseteqq\mathfrak{m}R^{\beta_{i_0-1}}$, with $i_0\geq 3$. Minimalizing, we get a resolution with $i_0$-th Betti number smaller than $\beta_{i_0}$. But from a bigraded free resolution of $\textrm{bigr}(N)$ we can deduce a graded free resolution of $\textrm{gr}^F(M)$, and the $i_0$-th Betti number of  $\textrm{gr}^F(M)$ is equal to the $i_0$-th Betti number of the statement, according to Proposition \ref{prop21}.
We can minimalize the resolution of $\textrm{gr}^F(M)$ obtained, which leads to a contradiction.

Now let us prove 2.
We have $\textrm{reg}_F(\textrm{ann}\delta/J)=1$. Indeed, the generator $t\tau+\chi$ has $F$-degree $1$, the other entries of the resolution we constructed have degree $0$ or $1$.
Similarly, we have $\textrm{reg}_F(R/J)=0$ because of the resolution we constructed.
Finally, in the proof of 1. we saw that the cone of the morphism $\alpha$ is the minimal resolution of $R/\textrm{ann}\delta$. Then
\[
\textrm{reg}_F N=\textrm{reg}_F\left(\frac{R}{\textrm{ann}\delta}\right)=
\textrm{max}\left(\textrm{reg}_F\left(\frac{\textrm{ann}\delta}{J}\right)-1,
\textrm{reg}_F\frac{R}{J}\right)=0.
\]

\subsection{Characterization of quasi-homogeneity}

Let $f\in\C\{x_1,\dots,x_n\}$ be an isolated singularity at $0$.
In this section, we will use the invariants coming from the minimal resolution of $N$ to characterize the quasi-homogeneity of $f$. 

First let us recall the decomposition of the module $N$ introduced in \cite{granger89}.
Let $E$ be a finite dimensional $\C$-vector space such that  $J(f)\oplus E=\C\{x\}$. Denote by $\D$ the ring of differential operators with analytic coefficients on $\C^{n}$, and  $\mathbb{D}$ the ring of differential operators with constant coefficients on $\C^{n}$, i.e. $\mathbb{D}=\sum_{\beta\in\mathbb{N}^n}\C\partial_x^{\beta}.$

\begin{prop}[\cite{granger89}, A.1.4]\label{prop17}
We have a decomposition of $N$ as a direct sum of $\C$-vector spaces
$$N=\D J(f)f^{s}\oplus(\oplus_{i\geq 0}\mathbb{D}E\partial_{t}^{i}f^{s})$$
and for any $i\geq 0$, the application
 $\mathbb{D}E\to \mathbb{D}E\partial_{t}^{i}f^{s}$
which maps $P\in\mathbb{D}E$ to $P\partial_{t}^{i}f^{s}$ is one-to-one.
\end{prop}

Define the left $\D$-ideal $I=\{P\in\D, Psf^s\in \D f^s\}$
and the $\mathcal{O}$-ideal 
 $\mathfrak{a}=\{u\in\mathcal{O},uf\in J(f)\}$.

\begin{lemme}\label{lemme4}
 $I=\D\mathfrak{a}$.
\end{lemme}
\paragraph{Proof}
Let $P=\sum \partial^{\beta}u_{\beta}$. 
We decompose
 $u_{\beta}f=a_{\beta}+b_{\beta}\in J(f)\oplus E.$
We have $sf^s=-f\partial_t f^s$ and $f'_i\partial_t f^s=-\partial_{x_i}f^s$, then
\[
 Psf^s=-\sum \partial^{\beta}u_{\beta}f\partial_t f^s=-\sum   \partial^{\beta}b_{\beta}\partial_t f^s\quad \textrm{mod}\, \D f^s.
\]
So $Psf^s\in\D f^s$ iff $\forall \beta, b_{\beta}=0$ after Proposition \ref{prop17}.$\square$

\begin{lemme}\label{lemme5}
 $F_{d-1}(\D)sf^s\cap \D f^s\subset F_d(\D)f^s$.
\end{lemme}
\paragraph{Proof}
Suppose $Psf^s=Qf^s$ with $P\in F_{d-1}(\D)$ and $d'=\textrm{ord}_F (Q)$ minimal. If $d'>d$, then 
$\sigma(Q)(f'_1,\dots,f'_n)=0.$
By Proposition \ref{prop16}, there exists $H\in F_{d'}(\D)$ such that $H f^s=0$ and $\sigma(H)=\sigma(Q)$. So $Qf^s=(Q-H)f^s$ and $\textrm{ord}_F(Q-H)<d'$, which contradicts the minimality of $d'$. $\square$

Let $u_1,\dots,u_l$ be a system of generators of $\mathfrak{a}$. There exist vector fields $\delta_i$ such that $\delta_i(f)=u_i f$.
We define then $Q_i=u_i s-\delta_i\in\textrm{ann}_{\D[s]}f^s$.
 
Let us recall the existence of a good operator anihilating $f^s$, i.e. an operator $S_L(s)\in\D[s]$, unitary and of degree $L$ in $s$, of the form $S_L(s)=\sum P_i s^i$ with $P_i\in F_{L-i}(\D)$ (cf. \cite{kashiwara}). We choose such an operator with $L$ minimal. 

\emph{Notation : } Let us denote for $i<j$, $S_{i,j}=f'_{i} \partial_{x_j}-f'_{j}\partial_{x_i}$. 

\begin{lemme}
 There exists an $F$-involutive system of generators of $\textrm{ann}_{\D[s]}f^s$ of the form $(S_{i,j})_{i<j},Q_1,\dots,Q_l,R_1,\dots,R_m,S_L$, where for all $i$, $R_i$ has degree in $s$ greater than $1$.
\end{lemme}
\paragraph{Proof}
By Proposition \ref{prop16}, the elements $S_{i,j}$ form an $F$-involutive system of generators of $\textrm{ann}_{\D}(f^s)$. Using this fact and Lemmas \ref{lemme4} and \ref{lemme5}, if $P(s)\in \textrm{ann}_{\D[s]}f^s$ has degree in $s$ less or equal to $1$, then 
\[
 P(s)\in \sum_{1\leq i\leq l} F_{\textrm{ord}_F (P(s))-1}(\D) Q_i +\sum_{i<j} F_{\textrm{ord}_F (P(s))-1}(\D)S_{i,j}.
\]
Take an $F$-involutive system of generators of $\textrm{ann}_{\D[s]}f^s$. We can replace the elements of degree in $s$ greater than $L$ by $S_L$, and those of degree in $s$ less than $2$ by $Q_1,\dots,Q_l,
(S_{i,j})$.$\square$

Now let us use $\textrm{Der(log $D$)}$, the set of vector fields $\chi$ such that $\chi (f)\in\mathcal{O}f$. Let us take a minimal system of generators $\delta_1,\dots,\delta_p$ of this $\mathcal{O}$-module. There exist functions $u_i$ such that $\delta_i (f)=u_i f$, so the $u_i$ generate $\mathfrak{a}$. Define  $V_i=u_i s-\delta_i\in\textrm{ann}_{\D[s]} f^s$.

\begin{lemme}
 There exists an $F$-involutive system of generators of $\textrm{ann}_{\D[s]}f^s$ of the form $V_1,\dots,V_p,R_1,\dots,R_m,S_L$, where for all $i$, $R_i$ has degree in $s$ greater than $1$.
\end{lemme}
\paragraph{Proof}
By the preceding lemma, it is sufficient to show that the $S_{i,j}$ are generated by the elements $V_i$, in a way adapted to the filtration $F$.
We have $S_{i,j}(f)=0$ so $S_{i,j}\in \textrm{Der(log $D$)}$, i.e. $S_{i,j}=\sum a_i \delta_i$.
Then $0=S_{i,j} (f)=\sum a_i \delta_i(f)=\sum a_i u_i f$ so $\sum a_i u_i=0$ and $S_{i,j}=-\sum a_i V_i$.$\square$

\emph{Notation : } Let us denote $W$ the ring $\textrm{gr}^V(\D_{x,t}^{(h)})$. 
For $P(s)\in \D[s]$, let $H(P(s))$ be its homogenization in $\mathbf{R}(\D[s])\simeq \textrm{gr}^V_0(\D_{x,t}^{(h)})\subset W$.

\begin{cor}\label{cor2}
 The set $$(f,(f'_i\partial_t),(H(V_i)),(H(R_i)),H(S_L))$$ is a bihomogeneous system of generators of $\textrm{ann}_W \delta$.
\end{cor}
It comes from the preceding lemma and Proposition \ref{prop4}.

\begin{lemme}\label{lemme14}
None of the elements $f,f'_i\partial_t,H(S_L)$ of this system can be removed.
\end{lemme}

\paragraph{Proof}
Recall that $W$ is a bigraded ring. Let $\bar{s}=-\partial_t t\in W$, and $k\geq 0$. The bigraded structure is defined by
\[
 W_{d,k}=\bigoplus_{0\leq i\leq d-k}\bar{s}^i (\D^{(h)})_{d-i-k}\partial_t^k
\quad\textrm{and}\quad
 W_{d,-k}=\bigoplus_{0\leq i\leq d}\bar{s}^i (\D^{(h)})_{d-i}t^k.
\]
Consider a bihomogeneous relation 
\begin{equation}\label{eq3.3}
 A f+\sum B_i f'_i\partial_t +\sum C_i V_i+\sum D_i H(R_i)+E H(S_L)=0.
\end{equation}
We have to show that none of the coefficients $A,B_i,E$ can be equal to $1$.
Considering the total degree, $A=1$ is impossible. 

Suppose $B_{i_0}=1$. Relation (\ref{eq3.3}) has then bidegree $(1,1)$. We must have for any $i,C_i=D_i=E=0$. It implies
\[
 f'_{i_0}\in \sum_{j\neq i_0}\mathcal{O}f'_j+\mathcal{O}f,
\]
which is impossible in the case of an isolated singularity $f$.

If $E=1$, considering the leading term in $\bar{s}$ in the relation (\ref{eq3.3}), we obtain that $1$ belongs to the maximal ideal of $W$, a contradiction.$\square$

\emph{Remark : } After the preceding lemma, we have $\beta_1>n+1.$ Since in the smooth case, $\beta_1=n+1$, we conclude that $\beta_1$ characterizes the smoothness (among germs for which $0$ is at most an isolated singularity).

We are now ready to characterize the quasi-homogeneity.
If $f$ is quasi homogeneous, the good operator has degree $1$ in $s$, is of the form $s-\theta$ with $\theta(f)=f$. From the resolution of $\textrm{bigr}N$ we constructed in Section 3.1, we deduce that   $\textrm{ann}_W \delta$ is minimally generated by $f,(f'_i\partial_t),(f'_i\partial_{x_j}-f'_j\partial_{x_i}),s-\theta.$ The corresponding shifts for the filtration $F$ are respectively $0,1,1,1$. 

If now $f$ is not quasi homogeneous, then $S_L$ has degree in $s$ greater than $1$. Indeed, suppose on the contrary that there exists a good operator anihilating $f^s$ of the form $s+P$ with $P=u+\theta$, $u\in\mathcal{O}$ and $\theta$ a vector field. We have
\[
 0=(s+P)f^s=(s(1+\theta(f)\frac{1}{f})+u)f^s.
\]
Substituting $s=0$, we get $u=0$. So $\theta(f)\frac{1}{f}=-1$ and $f\in J(f)$. This implies that $f$ is quasi homogeneous by \cite{saito71}.

\begin{prop}\label{prop20}
$f$ is quasi homogeneous if and only if $\textrm{reg}_F N=0$.
\end{prop}

\paragraph{Proof}
If $f$ is not quasi homogeneous, then $H(S_L)$ has degree in $F$ greater than $1$, but $H(S_L)$ is a necessary generator of $\textrm{ann}_W \delta$, so there exists $j$ such that $\mathbf{n}^{(1)}_j\geq 2$ and $\textrm{reg}_F N\geq 1$. 
If $f$ is quasi homogeneous, then $\textrm{reg}_F N=0$ by Theorem \ref{thm3}.$\square$

For plane curves, we can characterize the quasi-homogeneity by using $\beta_1$ only. 
If $f$ is quasi homogeneous, we have $\beta_1=5$. 

\begin{prop}\label{prop5}
 Suppose $n=2$ and $f$ not quasi homogeneous. Then $\beta_1\geq 6$.
\end{prop}
\paragraph{Proof}
We know that $\textrm{Der(log $D$)}$ is a free module of rank $2$ (cf. \cite{saito80}), so $p=2$.

It is sufficient to show that the following elements are necessary in the system of generators described in Corollary \ref{cor2} : $f, (f'_i\partial_t),(V_i),H(S_L)$. Take a bihomogeneous relation of type (\ref{eq3.3}).
By Lemma \ref{lemme14}, we only have to consider the case $C_i=1$.

Let us assume $C_1=1$. Considering the degrees, we have $D_i=E=0$, $C_2\in\mathcal{O}$, $B_i\in\mathcal{O}t$ and $A\in\mathcal{O}h\oplus\mathcal{O}t\partial_t\oplus \oplus \mathcal{O}\partial_{x_i}$. Let us take the terms in $\oplus \mathcal{O}\partial_{x_i}$ of  relation (\ref{eq3.3}) :
\begin{equation}
 \label{eq3.4}
 \delta_1\in \mathcal{O}\delta_2+\sum\mathcal{O}f\partial_i.
\end{equation}
Let us show that $f\partial_i\in\mathfrak{m}\textrm{Der(log $D$)}$.
Let
\[\left(
 \begin{array}{c}
\delta_1\\
\delta_2\\
\end{array}
\right)
=A \left(
\begin{array}{c}
\partial_1\\
\partial_2\\
\end{array}\right)
\]
with $A\in M_2(\mathcal{O})$. 
We have 
\[
 ^{t}(\textrm{co}A)
 \left(
 \begin{array}{c}
\delta_1\\
\delta_2\\
\end{array}
\right)
=\textrm{det}(A) \left(
\begin{array}{c}
\partial_1\\
\partial_2\\
\end{array}\right)
=u f \left(
\begin{array}{c}
\partial_1\\
\partial_2\\
\end{array}\right)
\]
with $u$ a unit, by \cite{saito80}. Moreover, for any $i,j, a_{i,j}\in \mathfrak{m}$. For if $a_{1,1}$ is a unit, 
then 
\[
 u_1 f=\delta_1(f)=\sum a_{1,j}f'_j
\]
so $f'_1\in\mathcal{O}f+\mathcal{O}f'_2$, which is impossible.

Coming back to (\ref{eq3.4}), we get 
$
 \delta_1\in\mathcal{O}\delta_2+\mathfrak{m}\textrm{Der(log $D$)},
$
which contradicts the minimality of the system $(\delta_i)$. $\square$

\end{document}